\documentclass[11pt]{article}
\usepackage{amsfonts}
\usepackage{amssymb}
\usepackage{amsmath}
=0pt
\def\sqr#1#2{{\vcenter{\vbox{\hrule height.#2pt
              \hbox{\vrule width.#2pt height#1pt \kern#1pt \vrule
width.#2pt}
              \hrule height.#2pt}}}}
\def\signed #1{{\unskip\nobreak\hfil\penalty50
              \hskip2em\hbox{}\nobreak\hfil#1
              \parfillskip=0pt \finalhyphendemerits=0 \par}}
\def\endpf{\signed {$\sqr69$}}
\def\dbE{{\mathbb{E}}}
\def\dbF{{\mathbb{F}}}

\def\dbP{{\mathbb{P}}}

\def\dbR{{\mathbb{R}}}

\def\si{\sigma}

\def\f{\varphi}

\def\3n{\negthinspace \negthinspace \negthinspace }
\def\2n{\negthinspace \negthinspace }
\def\1n{\negthinspace }
\def\ns{\noalign{\smallskip} }

\def\ds{\displaystyle}
\def\mmu{\mathcal{U}_{ad}}
\def\ups{\Upsilon_{\bar{u}}}
\def\ma{\mathcal{A}_{\bar{u}}}
\def\O{\Omega}
\def\Om{\Omega}
\def\cF{{\cal F}}

\def\cH{{\cal H}}

\def\cJ{{\cal J}}

\def\cL{{\cal L}}

\def\cQ{{\cal Q}}

\def\cU{{\cal U}}

\def\bu{\bar{u}}
\def\bx{\bar{x}}
\def\N{N^{b}}

\def\T{T^{b}}
\def\TT{T^{b(2)}}
\def\eps{\varepsilon}

\def\mr{\mathbb{R}}

\def\mmf{\mathbb{F}}

\def\mE{{\mathbb{E}}}

\def\no{\noindent}

\def\ms{\medskip}
\def\bs{\bigskip}
\def\q{\quad}
\def\qq{\qquad}
\def\hb{\hbox}

\def\limsup{\mathop{\overline{\rm lim}}}

\def\lan{\mathop{\langle}}
\def\ran{\mathop{\rangle}}

\def\wt{\widetilde}
\def\cd{\cdot}

\def\dim{\hbox{\rm dim$\,$}}

\def\ae{\hbox{\rm a.e.{ }}}

\def\({\Big (}
\def\){\Big )}
\def\[{\Big[}
\def\]{\Big]}

\def\={\buildrel \triangle \over =}

\def\resp{{\it resp. }}

\def\be{\begin{equation}}
\def\bel{\begin{equation}\label}
\def\ee{\end{equation}}
\def\bea{\begin{eqnarray}}
\def\eea{\end{eqnarray}}
\def\bt{\begin{theorem}}
\def\et{\end{theorem}}
\def\bc{\begin{corollary}}
\def\ec{\end{corollary}}
\def\bl{\begin{lemma}}
\def\el{\end{lemma}}
\def\bp{\begin{proposition}}
\def\ep{\end{proposition}}
\def\br{\begin{remark}}
\def\er{\end{remark}}
\def\ba{\begin{array}}
\def\ea{\end{array}}
\def\bd{\begin{definition}}
\def\ed{\end{definition}}

\newtheorem{lemma}{Lemma}[section]
\newtheorem{remark}{Remark}[section]

\newtheorem{theorem}{Theorem}[section]
\newtheorem{corollary}{Corollary}[section]

\newtheorem{definition}{Definition}[section]
\newtheorem{proposition}{Proposition}[section]

\makeatletter
   
   \@addtoreset{equation}{section}
\makeatother

\begin{document}

\title{\bf Necessary Conditions for Stochastic Optimal Control Problems in Infinite Dimensions}

\author{H\'el\`ene Frankowska\thanks{CNRS, IMJ-PRG, Sorbonne Universit\'e,
case 247, 4 place Jussieu,
75252 Paris, France.  The research of this
author partially benefitted from ``FMJH Program Gaspard
Monge in optimization and operation research", and from the support to this program from EDF
under grant  PGMO 2018-0047H, the  CNRS-NSFC PRC Project under grant 271392 and the AFOSR grant FA9550-18-1-0254. {\small\it E-mail:} {\small\tt helene.frankowska@imj-prg.fr}.\ms}~~~and~~~
Xu Zhang\thanks{School of Mathematics, Sichuan University, Chengdu 610064, China. The research of this
author is partially supported by NSF of China under grant 11221101, the NSFC-CNRS
Joint Research Project under grant 11711530142 and the PCSIRT under
grant IRT$\_$16R53 from the Chinese Education
Ministry. {\small\it E-mail:} {\small\tt
zhang$\_$xu@scu.edu.cn}.}}

\date{}

\maketitle

\begin{abstract}
The purpose of this paper is to establish the first and second order necessary conditions for stochastic optimal controls in infinite dimensions. The control system is governed by a stochastic evolution equation, in which both drift and diffusion terms may contain the control variable and the set of controls is allowed to be nonconvex. Only one adjoint equation is introduced to derive the first order optimality necessary condition either by means of the classical variational analysis approach or under some assumption which is quite natural in the deterministic setting to
guarantee the existence of optimal controls. More importantly, in order to avoid the essential difficulty  with the well-posedness of  higher order adjoint equations, using again the classical variational analysis approach, only the first and the second order adjoint equations are needed to formulate the second order necessary optimality condition, in which the solutions to the second order adjoint equation are understood in the sense of relaxed transposition.
\end{abstract}

\bs

\no{\bf 2010 Mathematics Subject
Classification}: 93E20, 49J53, 60H15

\bs

\no{\bf Key Words}: Stochastic optimal control; First and second order necessary optimality conditions;  Variational equation; Adjoint equation; Transposition solution; Maximum principle.


\section{Introduction and Notations}\label{s1}


Let $T>\tau\ge 0$, and
$(\O,\cF,\{\cF_t\}_{t\in [0,T]},\dbP)$ be a
complete filtered probability space (satisfying
the usual conditions). Write $\mathbf{F}=\{\cF_t\}_{t\in
[0,T]}$, and denote by $\dbF$ the progressive $\si$-field (in $[0,T]\times\Omega$) with respect to $\mathbf{F}$.
For  a Banach space $X$ with the norm $|\cdot |_{X}$, $X^*$ stands for the dual space of $X$.  For any
$t\in[0,T]$ and $r\in [1,\infty)$, denote by
$L_{\cF_t}^r(\O;X)$ the Banach space of all
$\cF_t$-measurable random variables $\xi:\O\to
X$ such that $\mathbb{E}|\xi|_X^r < \infty$,
with the canonical  norm.  Also, denote by
$D_{\dbF}([\tau,T];L^{r}(\O;X))$
the vector space of all $X$-valued, $r$-th power integrable
$\mathbf{F}$-adapted processes $\phi(\cd)$ such that $\phi(\cd):\; [\tau, T] \to L^r(\Omega,\cF_T , \dbP; X)$ is c\'adl\'ag, i.e.,
right continuous with left limits.
Then  $D_{\dbF}([\tau,T];L^{r}(\O;X))$ is a
Banach space with the norm
 $$
|\phi(\cd)|_{D_{\dbF}([\tau,T];L^{r}(\O;X))} \=
\sup_{t\in
[\tau,T)}\big(\mE |\phi(t)|_X^r\big)^{\frac{1}{r}}.
 $$
Similarly, denote by $C_{\dbF}([\tau,T];L^{r}(\O;X))$ the
Banach space of all $X$-valued $\mathbf{F}$-adapted
processes $\phi(\cdot)$ so that
$\phi(\cdot):[\tau,T] \to L^{r}_{\cF_T}(\O;X)$
is continuous, with the norm inherited from  $D_{\dbF}([\tau,T];L^{r}(\O;X))$.

Fix $r_1,r_2,r_3,r_4\in[1,\infty]$ and define
 $$
 \ba{ll}
\ds
L^{r_1}_\dbF(\O;L^{r_2}(\tau,T;X))=\Big\{\f:(\tau,T)\times\O\to
X\bigm|\f(\cd)\hb{
is $\mathbf{F}$-adapted and }\dbE\(\int_\tau^T|\f(t)|_X^{r_2}dt\)^{\frac{r_1}{r_2}}<\infty\Big\},\\
\ns\ds
L^{r_2}_\dbF(\tau,T;L^{r_1}(\O;X))=\Big\{\f:(\tau,T)\times\O\to
X\bigm|\f(\cd)\hb{ is $\mathbf{F}$-adapted and
}\int_\tau^T\(\dbE|\f(t)|_X^{r_1}\)^{\frac{r_2}
{r_1}}dt<\infty\Big\}.
 \ea
 $$
Both $L^{r_1}_\dbF(\O;L^{r_2}(\tau,T;X))$ and
$L^{r_2}_\dbF(\tau,T;L^{r_1}(\O;X))$ are Banach
spaces with the canonical norms. If $r_1=r_2$,
we simply denote the above spaces by
$L^{r_1}_\dbF(\tau,T;X)$.  For a Banach space $Y$, $\cL(X;Y)$ stands for  the (Banach)
space of all bounded linear operators from $X$
to $Y$, with the usual operator norm (When
$Y=X$, we simply write $\cL(X)$ instead of
$\cL(X, Y)$).

Let $H$ be an infinite dimensional separable Hilbert space, $V$ be another separable Hilbert space, and denote by $\cL_2^0$
the Hilbert space of all Hilbert-Schmidt
operators from $H$ to $V$. Let $W(\cdot)$ be a
$V$-valued, $\mathbf{F}$-adapted, standard $Q$-Brownian motion (with $Q\in \cL(V)$ being a positive definite, trace-class operator on $V$) or a cylindrical
Brownian motion. In the sequel, to simplify the presentation, we only consider the case of a
cylindrical Brownian motion.

Let $A$ be an unbounded linear operator
(with domain $D(A)\subset H$), which generates a
$C_0$-semigroup $\{S(t)\}_{t\geq 0}$ on $H$.
Denote by $A^*$ the dual operator of $A$. It is
well-known that $D(A)$ is a Hilbert space with
the usual graph norm, and $A^*$ generates a
$C_0$-semigroup $\{S^*(t)\}_{t\geq 0}$ on $H$,
which is the dual $C_0$-semigroup of
$\{S(t)\}_{t\geq 0}$.

Let $p\geq 1$ and denote by $p'$ the H\"older conjugate of $p$, i.e., $p'=\frac{p}{p-1}$ if $p\in (1,\infty)$; and $p'=\infty$ if $p=1$. Consider a nonempty, closed  subset $U$
of another separable Hilbert space $\wt H$. Put
 $$\cU_{ad}\triangleq \left\{u(\cdot)\in L_{\dbF}^p(\Omega; L^{2}(0, T;\widetilde H))\;\Big|\; u(t,\omega) \in U\hbox{ a.e. }(t,\omega)\in [0,T]\times \Omega\right\}.$$

Let us consider the following controlled
stochastic evolution equation:
\begin{equation}\label{ch-10-fsystem1}
\left\{
\begin{array}{lll}\ds
dx(t) =\big(Ax(t) +a(t,x(t),u(t))\big)dt + b(t,x(t),u(t))dW(t) \qq\mbox{in }(0,T],\\
\ns\ds x(0)=x_0,
\end{array}
\right.
\end{equation}
where $a(\cd,\cd,\cd):[0,T]\times H\times U \times \Omega\to H$ and
$b(\cd,\cd,\cd):[0,T]\times H\times U \times \Omega \to \cL_2^0$,   $u(\cd)\in \cU_{ad}$ and $x_0\in
H$. In \eqref{ch-10-fsystem1}, $u(\cd)$
is the control variable,
$x(\cd)=x(\cd;x_0,u(\cd))$ is the state
variable. As usual, when the context is clear, we omit the $\omega(\in\Omega)$ argument in the defined functions.

Let $f(\cd,\cd,\cd,\cd):[0,T]\times H\times U\times \Omega\to \dbR$. Define the cost functional $\cJ(\cdot)$ (for the
controlled system \eqref{ch-10-fsystem1}) as
follows:
\begin{equation}\label{jk1}
\cJ(u(\cdot))\triangleq \dbE\Big(\int_0^T
f(t,x(t),u(t))dt + g(x(T))\Big),\qquad\forall\,
u(\cdot)\in \cU_{ad},
\end{equation}
where $x(\cd)$ is the solution to
\eqref{ch-10-fsystem1}  corresponding to the control $u(\cdot)$.

We associate with  these data the following optimal control
problem for \eqref{ch-10-fsystem1}:

\ms

\no {\bf Problem (OP)} \q {\it Find a control $\bar
u(\cdot)\in \cU_{ad}$ such that
 \bel{jk2}
\cJ \big(\bar u(\cdot)\big) = \inf_{u(\cdot)\in
\cU_{ad}} \cJ \big(u(\cdot)\big).
 \ee
Any $\bar u(\cdot)$ satisfying (\ref{jk2}) is
called an {\it optimal control}. The
corresponding state process $\bar x(\cdot)$ is
called an {\it optimal state process}.
$\big(\bar x(\cdot),\bar u(\cdot)\big)$ is
called an {\it optimal pair}.}

The purpose of this paper is to establish first and second order necessary optimality conditions  for Problem (OP). We refer to  \cite{Bensoussan2, HP1, MM, TL, Zhou} and the references cited therein for some early works on the first order necessary optimality condition.
Nevertheless, all of these works on the necessary conditions for optimal
controls of infinite dimensional stochastic evolution equations addressed only some very special
cases, as  the case with the diffusion term not depending  on the control variable, or when the control domain is convex, and so on. For the general case when the diffusion term of the stochastic control system depends on  the control variable and the control region lacks convexity, stimulated by the work \cite{Peng90} addressing the same problems but in finite dimensions (i.e., $\dim H<\infty$), because presently $\dim H=\infty$, one has to handle a difficult problem   of the well-posedness of the second order adjoint equation, which is an operator-valued backward stochastic evolution equation (\ref{second ajoint equ}) (see the next section). This problem was solved at almost the same time in \cite{Du, Fuhrman, LZ2}, in which the Pontryagin-type maximum principles were established for optimal controls of general
infinite dimensional nonlinear stochastic systems. Nevertheless, the techniques used in \cite{LZ2} and \cite{Du, Fuhrman} are quite different. Indeed, since the most difficult part, i.e., the correction term ``$Q_2(t)$" in (\ref{second ajoint equ}) does not appear in the final formulation of the Pontryagin-type maximum principle, the strategy in both  \cite{Du} and  \cite{Fuhrman} is to ignore  completely $Q_2(t)$ in the well-posedness result for (\ref{second ajoint equ}). By contrast,  \cite{LZ2}  characterized the above mentioned $Q_2(t)$ because it was anticipated that this term should be useful somewhere else. More precisely, a new concept of solution, i.e., relaxed transposition solution was introduced in \cite{LZ2} to prove the well-posedness of (\ref{second ajoint equ}).

On the other hand, to the best of our knowledge, the only works on second order necessary optimality condition for Problem (OP) are \cite{lv2016, LZZ}, in which the control domains are supposed to be convex so that the classical convex variation technique can be employed. For the general stochastic optimal control problems in the finite dimensional framework,  when nonconvex control regions are considered and spike variations are
used as perturbations, as shown in \cite{ZZ2, ZZSIAMReview}, to derive the second order necessary optimality conditions, the cost functional needs to be expanded up to the forth order and four adjoint equations have to be introduced. Consequently,  an essential difficulty would arise in the present infinite dimensional situation. Indeed,  the infinite dimensional counterparts of multilinear function-valued backward stochastic differential equations introduced in \cite[(3.7)--(3.8)]{ZZ2} are   multilinear operator-valued backward stochastic evolution equations. So far the well-posedness of these sort of equations is completely unknown.

In this paper, in order to avoid the above mentioned difficulty with  the well-posednss, as in our earlier work dealing with the finite dimensional framework \cite{FZZ}, we use  the classical variational analysis approach to establish the
second order necessary optimality condition for Problem (OP) in the
general setting. The main advantage of this approach is that only two adjoint equations (i.e., (\ref{bsystem1}) and (\ref{second ajoint equ}) in the next section) are needed
to state the desired second order condition. It deserves mentioning that, thought $Q_2(t)$ in (\ref{second ajoint equ}) is useless in the first order optimality condition, it plays a key role in the statement of our second order necessary optimality condition (See the last two terms in (\ref{second order integral condition}) in Theorem \ref{TH second order integral condition}). Of course, as in \cite{FZZ}, the classical variational analysis approach can also be used to derive the
first order necessary optimality condition for Problem (OP). Furthermore, we also show a ``true" Pontryagin-type maximum principle for Problem (OP) under some mild assumption.

The rest of this paper is organized as follows. In Section \ref{s2}, some preliminaries on set-valued analysis and mild  stochastic processes are collected. Section \ref{s3}  is  devoted to the first order necessary optimality condition for Problem (OP); while in  Section \ref{s4}, we shall establish the second order necessary optimality condition. In Section \ref{s4.0} we prove the maximum principle under a convexity assumption using again a variational technique.

\section{Preliminaries}\label{s2}

\subsection{Some elements of set-valued analysis}

 For a nonempty subset $K\subset X$,
 the distance between a point  $x \in X$ and $K$ is defined  by $\displaystyle dist\,(x,K):= \inf_{y\in
K} |y-x|_{X}$. For any $\delta>0$, write $B(x,\delta)=\big\{y\in X\;\big|\;|y-x|_{X}<\delta\big\}$.

\begin{definition}
The adjacent cone $\T_{K}(x)$ to $K$ at $x$ is
defined by
$$\T_{K}(x):=\Big\{v\in X\ \Big|\ \lim_{\eps\to 0^+} \frac{dist\,(x+\eps v,K)}{\eps}=0         \Big\}.$$
\end{definition}

The dual cone of the tangent cone $\T_{K}(x)$, denoted by $\N_{K}(x)$, is called the normal cone of $K$ at $x$, i.e.,
$$\N_{K}(x):=\Big\{\xi\in X^*\ \Big|\ \xi(v)\le 0,\ \forall\; v\in \T_{K}(x) \Big\}.$$

\begin{definition}
For any $x\in K$ and $v\in \T_{K}(x)$, the second order adjacent subset to $K$ at $(x,v)$ is defined by
$$\TT_{K}(x,v):=\Big\{h\in X\ \Big|\ \lim_{\eps\to 0^+} \frac{dist\,(x+\eps v+\eps^{2}h,K)}{\eps^2}=0 \Big\}.$$
\end{definition}

Let $G:\ X \leadsto Y$ be  a
set-valued map with nonempty values. We denote by $graph (G)$  the
set of all $(x,y) \in X \times Y$ such that $y \in G(x)$.  Recall
that $G$ is called Lipschitz around $x$  if there exist $c \geq 0$
and $\delta
>0$  such that
$$ G(x_1) \subset G(x_2) + B(0,c|x_1-x_2|),\ \ \forall\;x_1, \; x_2 \in B(x,\delta).$$

 The
adjacent directional derivative $D^\flat G(x,y) (v)$  of $G$  at
$(x,y) \in graph (G)$ in a direction $v \in X$  is  a subset of
$Y$ defined by
$$w \in  D^\flat G(x,y) (v) \Leftrightarrow  (v,w) \in \T_{graph (G)}(x,y) .$$

If $G$  is locally Lipschitz at $x$, then
$$D^\flat G(x,y) (v) = \mbox{\rm Liminf}_{\varepsilon \to 0+}\frac{G(x+\varepsilon v) - y}{\varepsilon},$$
where $\mbox{\rm Liminf}$  stands for the Peano-Kuratowski limit (see for instance \cite{Aubin90}).

We shall need the following result.
\begin{proposition}{\rm (\cite[Proposition 5.2.6]{Aubin90})} \label{derivative}
Consider a set-valued map  $G :X \leadsto Y$ with nonempty convex
values and assume that $G$ is Lipschitz around $x$.  Then, for any
$y \in G(x)$, the values of  $D^\flat G(x,y)$  are convex and
$$D^\flat G(x,y)(0)= \T_{G(x)}(y)  \supset G(x) - y ,$$
$$D^\flat G(x,y)(v) + D^\flat G(x,y)(0)= D^\flat G(x,y)(v).$$
\end{proposition}

\subsection{Mild  Stochastic Processes}
In the rest of this paper, we shall denote by $C$ a generic
constant, depending on $T$, $A$, $p$ and $L$, which may be
different from one place to another.

First, we recall
that $x(\cd)\in C_\dbF([0,T];L^p(\O;H))$ is called a
mild solution to the equation \eqref{ch-10-fsystem1}
if
$$
x(t) = S(t)x_0 + \int_0^t
S(t-s)a(s,x(s),u(s))ds + \int_0^t
S(t-s)b(s,x(s),u(s))dW(s),\qq\forall\, t\in
[0,T].
$$

Throughout this section, we assume the following:

\ms

\no{\bf (A1)} {\it
$a(\cd,\cd,\cd):[0,T]\times H\times U \times \Omega\to H$ and
$b(\cd,\cd,\cd):[0,T]\times H\times U \times \Omega \to \cL_2^0$ are
two (vector-valued) functions such that

 {\rm i)} For
any $(x,u)\in H\times U$, the functions
$a(\cd,x,u,\cdot):[0,T] \times \Omega\to H$ and $b(\cd,x,u, \cd):[0,T] \times \Omega \to
\cL_2^0$ are $\dbF$-measurable;

{\rm ii)} For any
$(t,x)\in [0,T]\times H$, the functions
$a(t,x,\cd):U\to H$ and $b(t,x,\cd):U\to \cL_2^0$ are
continuous a.s.;

 {\rm iii)} There exist a constant $L>0$ and a nonnegative $\eta\in L^p(\Omega;L^1(0,T;\dbR))$
such that, for a.e. $(t,\omega)\in [0, T]\times\Omega$ and any $(x_1,x_2,u_1,u_2,u)\in
H\times H\times U\times U\times U$,
$$
\left\{
\begin{array}{ll}\ds
|a(t,x_1,u_1,\omega) - a(t,x_2,u_2,\omega)|_H+|b(t,x_1,u_1,\omega) -
b(t,x_2,u_2,\omega)|_{\cL_2^0} \\\ns\ds\q\leq
L(|x_1-x_2|_H+|u_1-u_2|_{\widetilde H}),  \\
\ns\ds |a(t,0,u,\omega)|_H +|b(t,0,u,\omega)|_{\cL_2^0}^2 \leq
\eta(t,\omega).
\end{array}
\right.
$$}

Further, we impose the following assumptions on the cost functional:

\ms

\no{\bf (A2)} {\it Suppose that
$f(\cd,\cd,\cd,\cd):[0,T]\times H\times U\times \Omega\to \dbR$
and $g(\cd,\cd):H\times \Omega\to \dbR$ are two functions
such that

{\rm i)} For any $(x,u)\in H\times U$, the
function $f(\cd,x,u,\cd):[0,T]\times \Omega\to \dbR$ is $\dbF$-measurable and $g(x,\cd):\Omega\to \dbR$ is $\cF_T$-measurable;

{\rm ii)} For any $(t,x)\in [0,T]\times
H$, the function $f(t,x,\cd):U\to \dbR$ is
continuous a.s.;

{\rm iii)} For all $(t,x_1,x_2,u)\in [0,T]\times H\times
H\times U$,
$$ \left\{
\begin{array}{ll}\ds
|f(t,x_1,u) - f(t,x_2,u)|_{H} +|g(x_1) -
g(x_2)|_H
 \leq L|x_1-x_2|_H,\ \ a.s.,\\
\ns\ds |f(t,0,u)|_H +|g(0)|_H \leq L,\ \ a.s.
\end{array}
\right.
$$}
It is easy to show the following result.

\begin{lemma}\label{estimatelinearsde}
Assume (A1). Then,  for any $x_{0}\in H$  and $u\in L_{\dbF}^p(\Omega; $ $L^{2}(0, T; \widetilde H))$, the state equation (\ref{ch-10-fsystem1}) admits a unique mild solution
$x \in C_\dbF([0, T];L^p(\Omega; H))$, and for any $t\in [0,T]$ the following estimate holds:
\begin{equation}\label{estimateofx}
\sup_{s\in[0,t]}\dbE |x(s,\cdot)|_H^p
\le C\dbE \Big[|x_{0}|_H^p
+\Big(\int_{0}^{t}|a(s,0,u(s),\cdot)|_Hds\Big)^p
+\Big(\int_{0}^{t}|b(s,0,u(s),\cdot)|_{\cL_2^0}^{2}ds\Big)^{\frac{p}{2}}\Big].
\end{equation}
Moreover, if  $\tilde{x}$ is the solution to (\ref{ch-10-fsystem1}) corresponding to
$(\tilde x_{0}, \tilde{u})\in H\times L_{\dbF}^p(\Omega; L^{2}(0, T;\widetilde H))$, then, for any $t\in [0,T]$,
\begin{equation}\label{estimateof delta x}
\sup_{s\in[0,t]}\dbE |x(s,\cdot)-\tilde{x}(s,\cdot)|_H^p
\le C\dbE ~\Big[|x_{0}-\tilde x_{0}|_H^p+\Big(\int_{0}^{t}|u(s,\cdot)
-\tilde{u}(s,\cdot)|_{\widetilde H}^{2}ds
\Big)^{\frac{p}{2}}\Big].
\end{equation}
\end{lemma}

For $\varphi=a,b, f$, denote
\bel{eq1205}
\varphi[t]=\varphi(t,\bar{x}(t),\bar{u}(t),\omega),\quad\varphi_{x}[t]=\varphi_{x}(t,\bar{x}(t),\bar{u}(t),\omega),\quad
\varphi_{u}[t]=\varphi_{u}(t,\bar{x}(t),\bar{u}(t),\omega)
\ee

and  consider the following
$H$-valued backward stochastic evolution
equation\footnote{Throughout this paper, for any
operator-valued process (\resp random variable)
$R$, we denote by $R^*$ its pointwise dual
operator-valued process (\resp random variable).
In particular, if $R\in L^1_\dbF(0,T; L^2(\Omega;
\cL(H)))$, then $R^*\in L^1_\dbF(0,T;
L^2(\Omega; \cL(H)))$, and $|R|_{L^1_\dbF(0,T;
L^2(\Omega; \cL(H)))}=|R^*|_{L^1_\dbF(0,T;
L^2(\Omega; \cL(H)))}$.}:
\begin{eqnarray}\label{bsystem1}
 \left\{
\begin{array}{l}
dP_{1}(t)=-  A^* P_1(t) dt-\big(a_{x}[t]^*P_{1}(t)
          +b_{x}[t]^*Q_{1}(t)
          -f_{x}[t]\big)dt+Q_{1}(t)dW(t), \quad  t\in[0,T], \\
P_{1}(T)=-g_{x}(\bar{x}(T)).
\end{array}\right.
\end{eqnarray}

Let us recall below the well-posedness result for the equation \eqref{bsystem1} in the
transposition sense, developed in \cite{LZ1, LZ2, LZ3}.

We consider the following (forward) stochastic
evolution equation:
\begin{equation}\label{fsystem1}
\left\{
\begin{array}{lll}\ds
dz = (Az + \psi_1(s))ds +  \psi_2 (s)dW(s) &\mbox{ in }(t,T],\\
\ns\ds z(t)=\eta,
\end{array}
\right.
\end{equation}
where $t\in[0,T]$,
$\psi_1\in L^1_{\dbF}(t,T;L^2(\O;H))$,
$\psi_2\in L^2_{\dbF}(t,T;L^2(\O; \cL_2^0))$ and
$\eta\in L^{2}_{\cF_t}(\O;H)$.
We call $(P_1(\cdot), Q_1(\cdot)) \in
D_{\dbF}([0,T];L^{2}(\O;H)) \times
L^2_{\dbF}(0,T;\cL_2^0))$  a transposition
solution to
 \eqref{bsystem1} if for any $t\in
[0,T]$, $\psi_1(\cdot)\in
L^1_{\dbF}(t,T;L^2(\O;H))$, $\psi_2(\cdot)\in
L^2_{\dbF}(t,T;L^2(\O; \cL_2^0))$, $\eta\in
L^2_{\cF_t}(\O;H)$ and the corresponding mild
solution $z\in C_{\dbF}([t,T];L^2(\O;H))$ to
\eqref{fsystem1}, the following is satisfied
\begin{equation}\label{eq def solzz}
\begin{array}{ll}\ds
\dbE \big\langle z(T),-g_x\big(\bar x(T)\big)\big\rangle_{H}
+\dbE\int_t^T \big\langle z(s),a_x[s]^* P_1(s) +b_x[s]^* Q_1(s)- f_x[s]\big\rangle_Hds\\
\ns\ds = \dbE \big\langle\eta,P_1(t)\big\rangle_H
+ \dbE\int_t^T \big\langle
\psi_1(s),P_1(s)\big\rangle_H ds + \dbE\int_t^T
\big\langle \psi_2(s),Q_1(s)\big\rangle_{\cL_2^0} ds.
\end{array}
\end{equation}

We have the
following result for the well-posedness of \eqref{bsystem1}.

\begin{lemma}\label{vw-th1}
{\rm (\cite{LZ2, LZ3})} The equation
\eqref{bsystem1} admits one and only one
transposition solution $(P_1(\cdot), $ $Q_1(\cdot))\in
D_{\dbF}([0,T];L^{2}(\O;H)) \times
L^2_{\dbF}(0,T;L^2(\O;\cL_2^0))$. Furthermore,
\begin{equation}\label{vw-th1-eq1}
|(P_1(\cdot), Q_1(\cdot))|_{
D_{\dbF}([0,T];L^{2}(\O;H)) \times L^2_{\dbF}(0,T;\cL_2^0))}\leq C\big(|g_x\big(\bar x(T)\big)|_{
L^2_{\cF_T}(\O;H)}+|f_x[\cd]|_{L^1_{\dbF}(0,T;L^2(\O;H))}\big).
\end{equation}
\end{lemma}

The proof of Lemma \ref{vw-th1} is based on the following Riesz-type Representation
Theorem established in \cite{LYZ}.

\begin{theorem}\label{rep}
Suppose $1\le q<\infty$, and that $X^*$ has
the Radon-Nikod\'ym property. Then
$$
\ds
L^p_\dbF(0,T;L^q(\Omega;X))^*=L^{p'}_\dbF(0,T;L^{q'}(\Omega;X^*)).
$$
\end{theorem}

Fix any
$r_1,r_2,r_3,r_4\in[1,\infty]$ and  denote by $
\cL_{pd}\big(L^{r_1}_{\dbF}(0,T;L^{r_2}(\Omega;X));L^{r_3}_{\dbF}(0,T;L^{r_4}(\Omega;Y))\big)$
the vector space of all bounded, pointwise
defined linear operators $\cL$ from
$L^{r_1}_{\dbF}(0,T;L^{r_2}(\Omega;X))$
to $L^{r_3}_{\dbF}(0,T;L^{r_4}(\Omega;Y))$, i.e.,
for $\ae (t,\omega)\in (0,T)\times\Omega$, there
exists an $L(t,\omega)\in\cL (X;Y)$ satisfying
$\big(\cL u(\cd)\big)(t,\omega)=L (t,\omega)u(t,\omega)$,
$\forall\, u(\cd)\in
L^{r_1}_{\dbF}(0,T;L^{r_2}(\Omega;X))$.

Consider the Hamiltonian
\begin{equation}\label{Hamiltonianconvex}
\mathbb{H} (t,x,u, v,w,\omega)
={\lan v,a(t,x,u,\omega)\ran}_H+{\lan w,b(t,x,u,\omega)\ran}_{\cL_2^0}-f(t,x,u,\omega),
\end{equation}
where $(t,x,u,v,w,\omega)\in [0,T]\times H\times \wt H\times H\times \cL_2^0\times\Omega$.
We then consider the following
$\cL(H)$-valued backward stochastic evolution
equation:
\begin{equation}\label{second ajoint equ}
\quad\left\{
\begin{array}{l}
dP_{2}(t)=-\Big(A^*P_{2}(t)+ P_{2}(t)A+a_{x}[t]^*P_{2}(t)+P_{2}(t)a_{x}[t] +b_{x}[t]^*P_{2}(t)b_{x}[t]\\
\qquad\qquad \qquad
 +b_{x}[t]^*Q_{2}(t)+Q_{2}(t)b_{x}[t]+\mathbb{H}_{xx}[t]\Big)dt+Q_{2}(t)dW(t),\qquad  t\in[0,T], \\
P_{2}(T)=-g_{xx}(\bar{x}(T)),
\end{array}\right.
\end{equation}
where $\mathbb{H}_{xx}[t]=\mathbb{H}_{xx}(t,\bar{x}(t),\bar{u}(t), P_{1}(t),Q_{1}(t))$ with $(P_1(\cdot),Q_1(\cdot))$ given by (\ref{bsystem1}).

Let us introduce
the following two (forward) stochastic evolution
equations:
\begin{equation}\label{op-fsystem1}
\left\{
\begin{array}{ll}
\ds dx_1 = (A+a_{x}[s])x_1ds + u_1(s)ds + b_{x}[s]x_1 dW(s) + v_1(s) dW(s) &\mbox{ in } (t,T],\\
\ns\ds x_1(t)=\xi_1
\end{array}
\right.
\end{equation}
and
\begin{equation}\label{op-fsystem2}
\left\{
\begin{array}{ll}
\ds dx_2 = (A+a_{x}[s])x_2ds + u_2(s)ds + b_{x}[s]x_2 dW(s) + v_2(s) dW(s) &\mbox{ in } (t,T],\\
\ns\ds x_2(t)=\xi_2,
\end{array}
\right.
\end{equation}
where
$
\xi_1,\xi_2 \in
L^{2p'}_{\cF_t}(\Omega;H)$, $u_1,u_2\in
L^2_\dbF(t,T;L^{2p'}(\Omega;H))$, and
$v_1,v_2\in
L^2_\dbF(t,T;$ $L^{2p'}(\Omega;\cL_2^0))$.
Also, we need to define the solution space
for \eqref{second ajoint equ}. For this purpose, for $p>1$, write
\begin{equation}\label{jshi1}
\begin{array}{ll}\ds
 L^{p}_{\dbF,w}(\Omega;D([0,T];\cL(H))\\
\ns\ds\= \Big\{P(\cd,\cd)\;\Big|\;
P(\cd,\cd)\in
\cL_{pd}\big(L_{\dbF}^{2}(0,T;L^{2p'}(\Omega;H));\;L^2_{\dbF}(0,T;L^{\frac{2p}{p+1}}(\Omega;H))\big),\\
\ns\ds\q P(\cd,\cd)\xi\!\in\!
D_{\dbF}([t,T];L^{\frac{2p}{p+1}}(\Omega;H))\;\mbox{ and }\;
\\
\ns\ds\q  |P(\cd,\cd)\xi|_{D_{\dbF}([t,T];L^{\frac{2p}{p+1}}(\Omega;H))}\leq C|\xi|_{L^{2p'}_{\cF_t}(\Omega;H)}\mbox{ for every } t\in[0,T]\hb{ and
}\xi\in L^{2p'}_{\cF_t}(\Omega;H) \Big\},
\end{array}
\end{equation}
 $$\cH_t\=L^{2p'}_{\cF_t}(\Om;H)\times
L^2_\dbF(t,T;L^{2p'}(\Om;H))\times
L^2_\dbF(t,T;L^{2p'}(\Om;\cL_2^0)),\q\forall\;t\in[0,T),
 $$
and\footnote{By Theorem \ref{rep}, since  $\cL_2^0$ is a Hilbert space, we deduce that $Q^{(t)}(0,0,\cd)^*$ is a bounded linear operator from $L^{2}_\dbF(t,T;$ $L^{\frac{2p}{p+1}}(\Om;\cL_2^0))^*=L^2_\dbF(t,T;L^{2p'}(\Om;\cL_2^0))$ to
$L^2_\dbF(t,T;L^{2p'}(\Om;\cL_2^0))^*=L^{2}_\dbF(t,T;L^{\frac{2p}{p+1}}(\Om;\cL_2^0))$. Hence, $Q^{(t)}(0,0,\cd)^*=\widehat
Q^{(t)}(0,0,\cd)$ makes sense.}
\begin{equation}\label{jshi2}
\begin{array}{ll}\ds \cQ^p[0,T]\ds \=\Big\{\big(Q^{(\cd)},\widehat
Q^{(\cd)}\big)\;\Big|\;
Q^{(t)},\widehat
Q^{(t)}\in\cL\big(\cH_t;\;
L^{2}_\dbF(t,T;L^{\frac{2p}{p+1}}(\Om;\cL_2^0))\big)\\
\ns\ds \qq\qq\q \mbox{and }Q^{(t)}(0,0,\cd)^*=\widehat
Q^{(t)}(0,0,\cd)\mbox{ for any } t\in [0,T)\Big\}.
\end{array}
\end{equation}

The notion of relaxed transposition solution
to \eqref{second ajoint equ} (\cite{LZ2, LZ3}):
We call
$
\big(P_2(\cd),\;Q_2^{(\cd)},$ $\widehat Q_2^{(\cd)}\big)\;\in\;
L^{p}_{\dbF,w}(\Omega;\;D([0,T];\; \cL(H)))\;\times$ $
\cQ^p[0,T]
$
a \index{relaxed transposition solution} relaxed transposition solution to the equation
\eqref{second ajoint equ} if for any $t\in [0,T]$,
$\xi_1,\xi_2\in L^{2p'}_{\cF_t}(\Omega;H)$,
$u_1(\cd), u_2(\cd)\in
L^2_{\dbF}(t,T;L^{2p'}(\Omega;H))$ and $v_1(\cd),
v_2(\cd)\in L^2_{\dbF}(t,T; L^{2p'}(\Omega;\cL_2^0))$,
the following is satisfied
\begin{equation}\label{6.18eq1}
\begin{array}{ll}
\ds\mE\big\langle -g_{xx}(\bar{x}(T)) x_1(T), x_2(T)
\big\rangle_{H} + \mE \int_t^T \big\langle
\mathbb{H}_{xx}[s] x_1(s), x_2(s) \big\rangle_{H}ds\\
\ns\ds =\mE\big\langle P_2(t) \xi_1,\xi_2
\big\rangle_{H}+ \mE \int_t^T
\big\langle P_2(s)u_1(s), x_2(s)\big\rangle_{H}ds
\\
\ns\ds \q+\mE\int_t^T\big\langle
P_2(s)x_1(s),
u_2(s)\big\rangle_{H}ds + \mE\int_t^T
\big\langle P_2(s)b_x(s)x_1 (s), v_2
(s)\big\rangle_{\cL_2^0}ds \\
\ns\ds \q+ \mE
\int_t^T\big\langle  P_2(s)v_1 (s), b_x (s)x_2 (s) + v_2(s)\big\rangle_{\cL_2^0}ds+ \mE\int_t^T
\big\langle v_1(s), \widehat
Q_2^{(t)}(\xi_2,u_2,v_2)(s)\big\rangle_{\cL_2^0}ds\\
\ns\ds \q+
\mE\int_t^T\big\langle
Q_2^{(t)}(\xi_1,u_1,v_1)(s), v_2(s)
\big\rangle_{\cL_2^0}ds,
\end{array}
\end{equation}
where $x_1(\cd)$ and $x_2(\cd)$ solve respectively
\eqref{op-fsystem1} and \eqref{op-fsystem2}.

We have the following well-posedness result for
the equation (\ref{second ajoint equ}) in the sense of relaxed transposition solution.

\begin{theorem}\label{OP-th2}
{\rm (\cite{LZ2, LZ3})} Assume that $p\in(1,2]$ and the Banach
space $L^p_{\cF_T}(\Omega;\dbR)$
is separable. Then, the
equation \eqref{second ajoint equ} admits one and
only one relaxed transposition solution $
\big(P_2(\cd),Q_2^{(\cd)},\widehat
Q_2^{(\cd)}\big)\in L^{p}_{\dbF,w}(\Omega;$ $D([0,T];
\cL(H)))\times \cQ^p[0,T]$. Furthermore,
\begin{equation}\label{ine the1z}
\begin{array}{ll}\ds
|P_2|_{ \cL(L_{\dbF}^{2}(0,T;L^{2p'}(\Omega;H));\;L^2(0,T;L_{\dbF}^{2p/(p+1)}(\Omega;H)))} + \sup_{t\in[0,T)}\big|\big(Q_2^{(t)},\widehat
Q_2^{(t)}\big)\big|_{\cL(\cH_t;\;
L^{2}_\dbF(t,T;L^{2p/(p+1)}(\Om;\cL_2^0)))^2}\\
\ns\ds\leq C\big(
|\mathbb{H}_{xx}[\cd])|_{L^1_\dbF(0,T;\;L^p(\Omega;\cL(H)))}
+
|g_{xx}(\bar x(T))|_{L^p_{\cF_T}(\Omega;\;\cL(H))}\big).
\end{array}
\end{equation}
\end{theorem}

\section{First order necessary optimality condition}\label{s3}

We impose the following further assumptions
on $a(\cd,\cd,\cd)$, $b(\cd,\cd,\cd)$,
$f(\cd,\cd,\cd)$ and $g(\cd)$.

\ms

\no{\bf (A3)} {\it The maps $a(t,\cd,\cd,\omega)$,
$b(t,\cd,\cd,\omega)$, and the functions $f(t,\cd,\cd,\omega)$ and
$g(\cd,\omega)$ are $C^1$ with respect to  $x$ and $u$.
Moreover, for a.e. $t\in [0,T]$ and any $(x,u)\in H\times
\wt H$,
\begin{equation}\label{ab}
\left\{
\begin{array}{ll}\ds
|a_x(t,x,u)|_{\cL(H)}+|b_x(t,x,u)|_{\cL(H;\cL_2^0)} + |f_x(t,x,u)|_H+|g_x(x)|_{H}\leq L,\ \ a.s.,\\
\ns\ds|a_u(t,x,u)|_{\cL(\wt H;H)}+
|b_u(t,x,u)|_{\cL(\wt H;\cL_2^0)} +
|f_u(t,x,u)|_{\wt H} \leq L,\ \ a.s.
\end{array}
\right.
\end{equation}}

\medskip

Now, let us introduce the classical first order variational control system. Let $\bu, v, v_{\eps}\in  L^{\beta}_{\mmf}(\Omega;$ $L^{2}(0,T;\wt H))$ ($\beta\ge 1$) satisfying $v_{\eps}\to v$ in $L^{\beta}_{\mmf}(\Omega;L^{2}(0,T;\wt H))$ as $\eps\to 0^+$. For $u^{\eps}:=\bu+\eps v_{\eps}$, let $x^{\eps}$ be the state  of (\ref{ch-10-fsystem1}) corresponding  to the control $u^{\eps}$, and put
$$\delta x^{\eps}=x^{\eps}-\bx.
$$

Consider the following linearized stochastic control system (Recall (\ref{eq1205}) for the notations $a_x[\cd]$, $a_u[\cd]$, $b_x[\cd]$ and $b_u[\cd]$):
\begin{equation}\label{first vari equ1}
\left\{
\begin{array}{l}
dy_{1}(t)= \big(Ay_{1}(t) +a_{x}[t] y_{1}(t)+a_{u}[t]v(t)\big)dt+\big(b_{x}[t] y_{1}(t)+ b_{u}[t]v(t)\big)dW(t),\quad  t\in[0,T], \\
y_{1}(0)=0.
\end{array}\right.
\end{equation}
The system (\ref{bsystem1}) is the first order adjoint equation associated with (\ref{first vari equ1}) and the cost function $f$. Similarly to \cite{FZZ}, it is easy to establish the following estimates.

\begin{lemma}\label{estimate one of varie qu}
If (A1) and (A3) hold, and  $\beta\ge 2$, then, for any $\bu, v,v_{\eps}$ and $\delta x^\eps$ as above
\begin{equation*}
 |y_{1}|_{L^\infty_{\dbF}(0,T;L^\beta(\O;H))}\le C |v|_{L^{\beta}_{\mmf}(\Omega;L^{2}(0,T;\wt H))},
\quad
|\delta x^{\eps}|_{L^\infty_{\dbF}(0,T;L^\beta(\O;H))}= O(\eps),\quad \hbox{as }\eps\to 0^+.
\end{equation*}
Furthermore,
\begin{equation}\label{r1 to 0}
|r_{1}^{\eps}|_{L^\infty_{\dbF}(0,T;L^\beta(\O;H))}\to 0,\quad \hbox{as }\eps\to 0^+,
\end{equation}
where
$$r_{1}^{\eps}(t,\omega):= \frac{\delta x^{\eps}(t,\omega)}{\eps}- y_{1}(t,\omega).
$$
\end{lemma}

Further, we have the following result.

\begin{theorem}\label{th max}
Let
the assumptions (A1) with $p=2$, (A2) and (A3) hold and
let $(\bar x(\cd),\bar u(\cd))$ be an optimal
pair of Problem (OP). Consider  the transposition solution  $(P_1(\cdot),Q_1(\cdot))$ to
\eqref{bsystem1}. Then,
\begin{equation}\label{first order integraltype condition}
\dbE\int_{0}^{T}{\lan \mathbb{H}_u[t],v(t)\ran}_{\wt H} dt\le 0,\qquad \forall\; v\in T^b_{{\cal U}_{ad}}(\bar u),
\end{equation}
where $\mathbb{H}_u[t]=\mathbb{H}_{u}(t,\bar x(t),\bar u(t),P_1(t),Q_1(t))$.
\end{theorem}

{\it Proof.}
Let $v\in \T_{{\cal U}_{ad}}(\bu(\cd))$. Then, for any $\eps>0$, there exists a  $v_{\eps}\in L_{\mmf}^{2}(\Omega;L^{2}(0,T;\wt H))$ such that $\bu+\eps v_{\eps}\in {\cal U}_{ad}$  and
$$\dbE\int^{T}_{0}|v(t)-v_{\eps}(t)|_{\wt H}^2dt\to 0,\ \hbox{ as }\eps\to 0^+.$$
Expanding the cost functional $J(\cdot)$ at $\bu(\cd)$, we have
\begin{eqnarray}\label{3.14}
0&\le& \frac{J(u^{\eps}(\cd))-J(\bu(\cd))}{\eps}\nonumber\\
&=& \dbE \int_{0}^{T}\Big(\int_{0}^{1}\lan f_{x}(t,\bx(t)+\theta\delta x^{\eps}(t),\bu(t)+\eps v_{\eps}(t)),\frac{\delta x^{\eps}(t)}{\eps}{\ran}_Hd\theta\nonumber\\
& &\qquad\qquad
+\int_{0}^{1}\lan f_{u}(t,\bx(t),\bu(t)+\theta \eps v_{\eps}(t)), v_{\eps}(t){\ran}_{\wt H}d\theta\Big)dt\nonumber\\
& &+\dbE \int_{0}^{1}\lan g_{x}(\bx(T)+\theta\delta x^{\eps}(T)),\frac{\delta x^{\eps}(T)}{\eps}{\ran}_Hd\theta\nonumber\\
&=&\dbE \int_{0}^{T}\big(\lan f_{x}[t],y_{1}(t){\ran}_H
+{\lan f_{u}[t],v(t)\ran}_{\wt H}\big)dt+\dbE{\lan g_{x}(\bx(T)),y_{1}(T)\ran}_H+\rho_{1}^{\eps},
\end{eqnarray}
where
\begin{eqnarray}\label{2016302e1}
\rho_{1}^{\eps}&=&\dbE \int_{0}^{T}\Big(\int_{0}^{1}\lan f_{x}(t,\bx(t)+\theta\delta x^{\eps}(t),\bu(t)+\eps v_{\eps}(t))
- f_{x}[t],\frac{\delta x^{\eps}(t)}{\eps}{\ran}_Hd\theta\nonumber\\
& &\qquad\quad
+\int_{0}^{1}{\lan f_{u}(t,\bx(t),\bu(t)+\theta \eps v_{\eps}(t))-f_{u}[t], v_{\eps}(t)\ran}_{\wt H}d\theta\nonumber\\
& &\qquad\quad
+\lan f_{x}[t],\frac{\delta x^{\eps}(t)}{\eps}-y_{1}(t){\ran}_H
+{\lan f_{u}[t],v_{\eps}(t)-v(t)\ran}_{\wt H}
\Big)dt\nonumber\\
& &+\dbE \int_{0}^{1}\lan g_{x}(\bx(T)+\theta\delta x^{\eps}(T))-g_{x}(\bx(T)),\frac{\delta x^{\eps}(T)}{\eps}{\ran}_Hd\theta\nonumber\\
& &\qquad\quad
+\dbE \lan g_{x}(\bx(T)),\frac{\delta x^{\eps}(T)}{\eps}-y_{1}(T){\ran}_H.
\end{eqnarray}
By Lemma \ref{estimate one of varie qu} (applied with $\beta =2$) and (A2)--(A3), using the Dominated
Convergence Theorem,  we conclude that
\begin{eqnarray*}
&&\Big|\dbE \int_{0}^{T}\int_{0}^{1}\lan f_{x}(t,\bx(t)+\theta\delta x^{\eps}(t),\bu(t)+\eps v_{\eps}(t))
- f_{x}[t],\frac{\delta x^{\eps}(t)}{\eps}{\ran}_Hd\theta dt\Big|\\
&&\le \Big(\dbE \int_{0}^{T}\int_{0}^{1}\big| f_{x}(t,\bx(t)+\theta\delta x^{\eps}(t),\bu(t)+\eps v_{\eps}(t))
- f_{x}[t]\big|_H^{2}d\theta dt\Big)^{\frac{1}{2}}\Big(\dbE \int_{0}^{T}\big|\frac{\delta x^{\eps}(t)}{\eps}  \big|_H^{2}dt\Big)^{\frac{1}{2}}\\
& &\quad\to 0,\quad \hbox{ as }\eps\to 0^+.
\end{eqnarray*}
Similarly, we have
\begin{eqnarray*}
\dbE \int_{0}^{T}\int_{0}^{1}{\lan f_{u}(t,\bx(t),\bu(t)+\theta \eps v_{\eps}(t))-f_{u}[t], v_{\eps}(t)\ran}_{\wt H}d\theta dt\to 0,\quad \hbox{ as } \eps\to 0^+,
\end{eqnarray*}
and
\begin{eqnarray*}
\dbE \int_{0}^{1}\lan g_{x}(\bx(T)+\theta\delta x^{\eps}(T))-g_{x}(\bx(T)),\frac{\delta x^{\eps}(T)}{\eps}{\ran}_Hd\theta\to 0,\quad \hbox{ as } \eps\to 0^+.
\end{eqnarray*}
Then, by (A2)--(A3) and  Lemma \ref{estimate one of varie qu},
we obtain that
\begin{eqnarray}\label{2016302e2}
\lim_{\eps\to 0^+}\big|\rho_{1}^{\eps}\big|
&\le&\limsup_{\eps\to 0^+}\Big|\dbE \int_{0}^{T}\lan f_{x}[t],\frac{\delta x^{\eps}(t)}{\eps}-y_{1}(t){\ran}_Hdt\Big|\nonumber\\
& &
+\limsup_{\eps\to 0^+}\Big|\dbE \int_{0}^{T}{\lan f_{u}[t],v_{\eps}(t)-v(t)\ran}_{\wt H}
dt\Big|\nonumber\\
& &+\limsup_{\eps\to 0^+}\Big|\dbE \lan g_{x}(\bx(T)),\frac{\delta x^{\eps}(T)}{\eps}-y_{1}(T){\ran}_H\Big|= 0.
\end{eqnarray}
Therefore, from (\ref{3.14}) and (\ref{2016302e2}), we conclude that
\begin{eqnarray}\label{first order taylor exp}
0&\le& \dbE \int_{0}^{T}\big({\lan f_{x}[t],y_{1}(t)\ran}_H
+{\lan f_{u}[t],v(t)\ran}_{\wt H}\big)dt +\dbE{\lan g_{x}(\bx(T)),y_{1}(T)\ran}_H.
\end{eqnarray}

By means of the definition of transposition
solution to
 \eqref{bsystem1}, we have
\begin{equation}\label{duality between y1 p1}
\ba{ll}
\ds-\dbE {\lan g_{x}(\bar{x}(T)),y_{1}(T)\ran}_H
=\dbE{\lan P_{1}(T),y_{1}(T)\ran}_H\\[3mm]\ds
=\dbE\int_{0}^{T}\big({\lan P_{1}(t),a_{x}[t]y_{1}(t)\ran}_H
+{\lan P_{1}(t),a_{u}[t]v(t)\ran}_H\\[3mm]\ds
\qquad\quad+{\lan Q_{1}(t),b_{x}[t]y_{1}(t)\ran}_{\cL_2^0}
+{\lan Q_{1}(t),b_{u}[t]v(t)\ran}_{\cL_2^0}\\[3mm]\ds
\ds\qquad\quad- {\lan a_{x}[t]^*P_{1}(t),y_{1}(t)\ran}_H
-{\lan b_{x}[t]^*Q_{1}(t),y_{1}(t)\ran}_H
+{\lan f_{x}[t],y_{1}(t)\ran}_H\big)dt\\[3mm]\ds
\ds=\dbE\int_{0}^{T}\big({\lan P_{1}(t),a_{u}[t]v(t)\ran}_H
+{\lan Q_{1}(t),b_{u}[t]v(t)\ran}_{\cL_2^0}+{\lan f_{x}[t],y_{1}(t)\ran}_H
\big)dt.
\ea
\end{equation}

Substituting (\ref{duality between y1 p1}) in (\ref{first order taylor exp}), and recalling \eqref{Hamiltonianconvex}, we obtain that
\begin{eqnarray}\label{first order vari inequ}
0&\le&  -\dbE \int_{0}^{T}\big({\lan P_{1}(t),a_{u}[t]v(t)\ran}_H
+{\lan Q_{1}(t),b_{u}[t]v(t)\ran}_{\cL_2^0}
-{\lan f_{u}[t],v(t)\ran}_{\wt H}\big)dt\nonumber\\
&=&-\dbE\int_{0}^{T}{\lan \mathbb{H}_{u}[t],v(t)\ran}_{\wt H}dt,
\end{eqnarray}
which gives (\ref{first order integraltype condition}).
This completes the proof of Theorem \ref{th max}.
\endpf

\ms

As a consequence of Theorem \ref{th max}, one has  the following pointwise first order necessary condition for optimal pairs of Problem (OP) that follows easily from measurable selection theorems.

\begin{theorem}\label{TH first order pointwise condition}
Under the assumptions in Theorem \ref{th max}, it holds that
\begin{equation}\label{first order pointwise condition}
 \mathbb{H}_{u}[t] \in \N_{U}(\bu(t)),\ a.e.\ t\in [0,T],\ a.s.
\end{equation}
\end{theorem}

\section{Maximum Principle}\label{s4.0}

In this section, we address the Pontryagin-type maximum principle for Problem (OP). For this purpose, we
introduce the following assumption.

\ms

\no{\bf (A4)} {\it For a.e. $(t,\omega) \in
[0,T] \times \Omega$  and  for all $x \in H$, the set
$$F(t,x,\omega):=\big\{(a(t,x,u,\omega), b(t,x,u,\omega), f(t,x,u,\omega) +r)\;\big|\; u \in U, \; r \geq 0\big\}$$
is closed and convex in $H\times H\times \dbR$.}

\ms

The above  assumption is familiar in the deterministic optimal control
where it  is very useful  to guarantee the existence of optimal
controls. In particular, it holds true whenever $a,b$ are affine in the control, $U$ is convex and compact and $f$ is convex with respect to $u.$

Also, we impose  the following assumption, which is  weaker than (A1) because Lipschitz continuity with respect to $u$ is no longer required.

\ms

\no{\bf (A5)} {\it Suppose that $a:[0,T]\times
H\times U \times \Omega \to H$ and $b:[0,T]\times H\times U \times \Omega \to \cL_2^0$
are two (vector-valued) functions satisfying the conditions {\rm i)} and {\rm ii)} in {\rm (A1)}, and such that for some $L \geq 0,$  $\eta \in L^2([0,T] \times \Omega; \dbR)$ and  for a.e.
$(t,\omega)\in [0, T]\times\Omega$ and any $(x_1,x_2,u)\in
H\times H\times U$,
$$
\left\{
\begin{array}{ll}\ds
|a(t,x_1,u,\omega) - a(t,x_2,u,\omega)|_H+|b(t,x_1,u,\omega) -
b(t,x_2,u,\omega)|_{\cL_2^0} \leq
L|x_1-x_2|_H,\\
\ns\ds |a(t,0,u,\omega)|_H +|b(t,0,u,\omega)|_{\cL_2^0}^2 \leq
\eta(t,\omega).
\end{array}
\right.
$$}

Under the assumption (A5), the system
\eqref{ch-10-fsystem1} is  still well-posed.

We have the following Pontryagin-type maximum principle for Problem (OP) (Recall (\ref{Hamiltonianconvex}) for the definition of $\mathbb{H}(\cd)$):

\begin{theorem} \label{PMPconvex} Assume  (A2), (A4) and (A5).  If  $(\bar x, \bar u)$ is an optimal pair for Problem (OP), then $(P_1,Q_1)$  defined in Lemma \ref{vw-th1}  verifies for a.e. $t \in [0,T]$ the maximality condition
$$ \max_{u \in U} \mathbb{H}(t,\bar x(t),u, P_{1}(t),Q_{1}(t))=\mathbb{H}(t,\bar x(t),\bar u(t), P_{1}(t),Q_{1}(t)),\;\; a.s.$$
\end{theorem}
{\it Proof.}  Fix $u (\cd)\in {\mathcal U}_{ad}$
and  consider the following linearized system:
\begin{equation}\label{first vari equ}
\left\{
\begin{array}{l}
d\tilde y_{1}(t)= \big(A\tilde y_{1}(t) +a_{x}[t] \tilde y_{1}(t)+ \alpha(t)\big)dt+\big(b_{x}[t] \tilde y_{1}(t)+ \beta(t)\big)dW(t),\quad  t\in[0,T], \\
\tilde y_{1}(0)=0,
\end{array}\right.
\end{equation}
where
$$\alpha (t) = a(t,\bar x(t),u(t)) - a [t],\qquad \beta (t) = b(t,\bar x(t),u(t)) - b [t].$$
Then
$$|\alpha (t,\omega)|_H + |\beta (t,\omega)|_{\cL_2^0} \leq C(|\bar x(t,\omega)| + |\eta (t,\omega)|+1 ).$$

Denote by $D_x^\flat F$ the adjacent derivative of $F$ with respect to $x$.
One can easily check  that for a.e. $t \in [0,T]$,
$$(a_{x}[t]\tilde y_1(t),b_{x}[t]\tilde y_1(t), f_{x}[t]\tilde y_1(t) ) \in D^\flat _x F(t,\bar x(t), a [t],b [t],f [t]   )(\tilde y_1(t)), \;\;  a.s. $$

Define
$$\gamma (t) = f(t,\bar x(t),u(t)) - f [t].$$
By Proposition \ref{derivative},  we have
$$ \big(a_{x}[t] \tilde y_{1}(t)+ \alpha(t),b_{x}[t] \tilde y_{1}(t)+ \beta(t), f_{x}[t] \tilde y_{1}(t)+ \gamma(t) \big)
\in D_x^\flat F(t,\bar x(t), a [t],b [t],f
[t]   )(\tilde y_1(t)).$$

Recall that $\dbF$ stands for the progressive $\si$-field (in $[0,T]\times\Omega$) with respect to $\mathbf{F}$. Let us consider the measure space $([0,T]\times\Omega,  \dbF, dt\times d \dbP)$ and its completion $([0, T]\times \Omega, {\mathcal A}, \mu )$.

Fix any sequence  $\varepsilon_i \downarrow 0$  and define for every $j \geq 1$ the sets
$$
\ba{ll}\widehat H^i_j(t,\omega)= \big\{(u,r) \in U \times \dbR_+\;\big|\;\\ \qq\qq\q
|a(t,\bar x(t,\omega) +\varepsilon_i \tilde y_1(t,\omega) ,u,\omega) - (a[t]
 +\varepsilon_i(a_x[t]\tilde y_1(t,\omega)+\alpha (t,\omega)))|_H  \leq \varepsilon_i /2^j \\ \qq\q\qq  |b(t,\bar x(t,\omega) +\varepsilon_i \tilde y_1[t,\omega] ,u,\omega)
  -(b[t]+ \varepsilon_i(b_x[t]\tilde y_1(t,\omega)+\beta (t,\omega)))|_{\cL_2^0} \leq \varepsilon_i /2^j,  \\ \qq\q\qq |f(t,\bar x(t,\omega) +\varepsilon_i \tilde y_1(t,\omega) ,u,\omega) +r - (f[t]  +\varepsilon_i(f_x[t]\tilde y_1(t,\omega)+\gamma (t,\omega)))|  \leq \varepsilon_i /2^j  \big\}
  \ea
$$
and
$$H^i_j (t,\omega)= \left\{ \begin{array}{ll} \widehat H^i_j (t,\omega) & {\rm if } \;\; \widehat H^i_j (t,\omega) \neq \emptyset \\
\{(\bar u(t,\omega),0)\} \quad& {\rm otherwise. }\end{array} \right.$$

By the definition of set-valued derivative, we can find  a subsequence ${i_j}$ and  a decreasing family of measurable sets $A_j \subset [0,T] \times \Omega$, such that
 $\lim_{j \to \infty}\mu (A_j)=0$ and   $\widehat H^{i_j}_j (t,\omega)\neq \emptyset$  for $(t,\omega) \in ([0,T] \times \Omega) \, \backslash \,   A_j$.

By \cite[Theorem 8.2.9]{Aubin90},   $ H^{i_j}_j$  is
$\mu$-measurable.  Since it also has closed nonempty values,  it
admits a measurable selection.  Modifying this selection on a set
of measure zero, we obtain $u_j \in {\mathcal U}_{ad}$ and an $\mathbf{F}$-adapted, real-valued process $r_j$  satisfying
 for a.e.  $(t,\omega)$ with $\widehat
H^{i_j}_j (t,\omega)\neq \emptyset$  the following inequalities:
$$
\ba{ll}
\ds
|a(t,\bar x(t,\omega) +\varepsilon_{i_j} \tilde y_1(t,\omega) ,u_j(t,\omega),\omega) - (a[t]  +\varepsilon_{i_j}(a_x[t]\tilde y_1(t,\omega)+\alpha (t,\omega)))|_H
 \leq \varepsilon_{i_j}/2^j,\\
 \ns\ds  |b(t,\bar x(t,\omega) +\varepsilon_{i_j} \tilde y_1(t,\omega) ,u_j(t,\omega),\omega) -(b[t]  +
 \varepsilon_{i_j}(b_x[t]\tilde y_1(t,\omega)+\beta (t,\omega)))|_{\cL_2^0}  \leq \varepsilon_{i_j}/2^j,  \\
  \ns\ds |f(t,\bar x(t,\omega) +\varepsilon_{i_j} \tilde y_1(t,\omega) ,u_j(t,\omega),\omega)  +r_j(t,\omega) - (f[t] +\varepsilon_{i_j}(f_x[t]\tilde y_1(t,\omega)+\gamma (t,\omega)))|  \leq \varepsilon_{i_j}/2^j.
  \ea
$$

Define
$$
\ba{ll}
\ds a_j (t) = a(t,\bar x(t)+\varepsilon_{i_j}\tilde y_1(t), u_j (t)), \qquad b_j (t)= b(t,\bar
x(t)+\varepsilon_{i_j}\tilde y_1(t), u_j (t)), \\
\ns\ds f_j(t)=f(t,\bar x(t)+\varepsilon_{i_j}\tilde y_1(t), u_j (t)).
\ea
$$

Observe that if  $\widehat
H^{i_j}_j (t,\omega)=\emptyset$, then by (A2), (A5), it follows that
$$
\ba{ll}
\ds |a_j(t,\omega)  - (a[t,\omega]  +\varepsilon_{i_j}(a_x[t,\omega]\tilde y_1(t,\omega)+\alpha (t,\omega)))|_H
 \leq C \varepsilon_{i_j}(|\tilde y_1(t,\omega) | + |\bar x(t,\omega)| +|\eta(t,\omega)|), \\
 \ns\ds
 |b_j(t,\omega)  - (b[t,\omega]  +\varepsilon_{i_j}(b_x[t,\omega]\tilde y_1(t,\omega)+\beta (t,\omega)))|_{\cL_2^0}
 \leq C \varepsilon_{i_j}(|\tilde y_1(t,\omega) | + |\bar x(t,\omega)| +|\eta(t,\omega)|), \\
 \ns\ds |f_j(t,\omega)  - (f[t,\omega]  +\varepsilon_{i_j}(f_x[t,\omega]\tilde y_1(t,\omega)+\gamma(t,\omega)))|
 \leq C \varepsilon_{i_j}(|\tilde y_1(t,\omega) | + |\bar x(t,\omega)| + L).
 \ea
 $$

Consider the solution $x_j$ of the following stochastic equation:
\begin{equation}\label{ch-10-fsystem1j0}
\left\{
\begin{array}{lll}\ds
dx_j(t) =\big(Ax_j(t) +a(t,x_j(t),u_j(t))\big)dt + b(t,x_j(t),u_j(t))dW(t) \qq\mbox{in }(0,T],\\
\ns\ds x_j(0)=x_0.
\end{array}
\right.
\end{equation}
 Then, by (\ref{ch-10-fsystem1j0}) we have
\begin{equation}\label{ch-10-fsystem1j}
\begin{array}{l}
\ds\dbE |x_j(t)- \bar x(t) - \varepsilon_{i_j} \tilde y_1(t)|_H^2 \\[3mm]
\ds\leq 2\dbE
\left|\int_0^t S(t-s)(a(s,x_j(s),u_j(s)) - a[s] -
\varepsilon_{i_j}(a_x[s]\tilde y_1(s)+\alpha (s)) )ds\right|_H^2 \\[3mm]
\q\ds+  2 \dbE \left|\int_0^t S(t-s)(b(s,x_j(s),u_j(s)) - b[s] -
\varepsilon_{i_j}(b_x[s]\tilde y_1(s)+\beta (s)))dW(s)\right|_H^2 \\[3mm]
\ds\leq C \dbE \int_0^t |a_j(s) - a[s] -
\varepsilon_{i_j}(a_x[s]\tilde y_1(s)+\alpha (s)) |_H^2ds \\[3mm]
\q \ds+ C \dbE
\int_0^t
|x_j(s)- \bar x(s) - \varepsilon_{i_j}\tilde y_1(s) |_H^2ds \\[3mm]
\ds\q+ C \dbE \int_0^t |b(s,x_j(s),u_j(s)) - b[s] -
\varepsilon_{i_j}(b_x[s]\tilde y_1(s)+\beta (s))|_{\cL_2^0}^2ds)
\\[3mm]\ds
\leq C \varepsilon_{i_j}^2 2^{-j} + C\varepsilon_{i_j}^2 \int_{A_j}\left(
|\tilde y_1(s)|_H^{2} +|\bar x(s)|_H^{2} + |\eta (s)|^{2} \right)ds \, d\omega  \\[3mm]
\ds\q  + C  \int_0^t \dbE
|x_j(s)- \bar x(s) - \varepsilon_{i_j} \tilde y_1(s)|_H^2ds.
\end{array}
\end{equation}
This and the Gronwall inequality imply that
for a constant $C>0$, a sequence $\delta_j \downarrow 0$  and  for all $j$

$$ \dbE |x_j(t)- \bar x(t) - \varepsilon_{i_j} \tilde y_1(t)|_H^2 \leq C\varepsilon_{i_j}^2 (2^{-j} +\delta_j). $$
Consequently

$$\lim_{j \to \infty} \sup_{t \in [0,T]}  \dbE  \left|\frac{x_j(t) - \bar x(t)}{\varepsilon_{i_j}} - \tilde y_1(t)\right|_H^2 =  0 . $$

On the other hand
\begin{equation}\label{ch-10-fsystem1f}
\begin{array}{l}
\ds\dbE \left|\int_0^T (f(s,x_j(s),u_j(s)) + r_{\varepsilon_{i_j}}(s) -
f[s]
- \varepsilon_{i_j}f_x[s]\tilde y_1(s) )ds\right|  \\[3mm]\ds
\leq  \dbE \int_0^T \left|f_j(s) + r_{\varepsilon_{i_j}}(s) - f[s] - \varepsilon_{i_j}f_x[s]\tilde y_1(s) \right|ds + C\dbE \int_0^T
\left|\bar x(s) +\varepsilon_{i_j} \tilde y_1(s)  - x_j(s) \right|_Hds .
\end{array}
\end{equation}

Thus
$$\lim_{j \to \infty} \frac{1}{\varepsilon_{i_j}} \dbE \left|\int_0^T (f(s,x_j(s),u_j(s)) + r_{\varepsilon_{i_j}}(s) - f[s]
- \varepsilon_{i_j}f_x[s]\tilde y_1(s) )ds\right| =0. $$

By the optimality of $(\bar x, \bar u)$, we arrive at
\begin{equation}\label{ch-10-fsystem1comp}
\begin{array}{l}
\ds 0 \leq \dbE \left[\int_0^T \left(f(s,x_j(s),u_j(s)) +
r_{\varepsilon_{i_j}}(s) - f[s]\right)ds + g(x_j(T)) - g(\bar x(T))\right] \\[3mm]\ds
\q\leq \varepsilon_{i_j}\int_0^T {\lan f_x[s],\tilde y_1(s) \ran}_Hds \\[3mm]\ds\qq+ \dbE \left|\int_0^T
\left(f(s,x_j(s),u_j(s)) + r_{\varepsilon_{i_j}}(s) - f[s] -
\varepsilon_{i_j}{\lan f_x[s],\tilde y_1(s) \ran}_H\right)ds\right|
 + g(x_j(T)) - g(\bar x(T)).
\end{array}
\end{equation}
Dividing by $\varepsilon_{i_j}$  and taking the limit yields
\begin{eqnarray}
0&\le& \dbE \int_{0}^{T}{\lan f_{x}[t],\tilde y_{1}(t)\ran}_H dt
+\dbE{\lan g_{x}(\bx(T)),\tilde y_{1}(T)\ran}_H.
\end{eqnarray}

By means of the definition of transposition
solution to
 \eqref{bsystem1}, in the same way as that in (\ref{duality between y1 p1})  we
finally obtain that
 \begin{equation}\label{integineq}
 \ba{ll}\ds
\dbE \int_{0}^{T}\big({\lan P_{1}(t),a(t,\bar x(t),u(t)) - a
[t]\ran}_H +{\lan Q_{1}(t),b(t,\bar x(t),u(t)) - b
[t]\ran}_{\cL_2^0}\\
\ns\ds\qq\q -(f(t,\bar x(t),u(t)) - f [t])\big)dt \\
\ns\ds\leq 0.
\ea
 \end{equation}

For  all $(t,\omega) \in [0,T] \times \Omega$ and integer $i \geq 1$ define the sets
$$
\ba{ll}\ds
K_i (t,\omega)  = \big\{u \in U\;\big|\; {\lan P_{1}(t,\omega),a(t,\bar x(t,\omega),u,\omega) - a
[t]\ran}_H +{\lan Q_{1}(t,\omega),b(t,\bar x(t,\omega),u,\omega) - b
[t]\ran}_{\cL_2^0}\\[3mm]
\ds\qq\qq\qq \qq\q-(f(t,\bar x(t,\omega),u,\omega) - f [t]) \geq  2^{-i} \big\}
\ea
$$
and
$$\widehat K_i(t,\omega) = \left\{ \begin{array}{ll}   K_i(t,\omega) \qq& {\rm if }\;\; K_i(t,\omega) \neq \emptyset, \\ \{\bar u (t,\omega)\} & {\rm otherwise}. \end{array}\right.$$
Let $A_i$ denote the set of all $(t,\omega)$  such that $K_i(t,\omega)\neq \emptyset$.
To end the proof it is enough to show that $\mu-$measure of  the following set
$$
\ba{ll}
\ds\big\{(t,\omega) \;\big|\; \exists \; u \in U \;\mbox{\rm so that}\;  {\lan P_{1}(t,\omega),a(t,\bar x(t,\omega),u,\omega) - a
[t]\ran}_H  \\
\ns\ds\qq\qq\qq\qq\qq\;\,\; +{\lan Q_{1}(t,\omega),b(t,\bar x(t,\omega),u,\omega) - b
[t]\ran}_{\cL_2^0} -(f(t,\bar x(t,\omega),u,\omega) - f [t] > 0\big\}
\ea $$
is equal to zero.
 Indeed otherwise there would exist an $i \geq 1$ such that $\mu(A_i) >0.$
   By \cite[Theorem 8.2.9]{Aubin90}  $ K_i(\cdot,\cdot)$ admits a $\mu-$measurable selection.  Modifying it on a set of measure zero, we obtain a  control  $u (\cd)\in {\mathcal U}_{ad}$  such that for a.e. $(t,\omega) \in A_i$,
 $$
 \ba{ll}\ds{\lan P_{1}(t,\omega),a(t,\bar x(t,\omega),u(t,\omega),\omega) - a
[t]\ran}_H +{\lan Q_{1}(t,\omega),b(t,\bar x(t,\omega),u(t,\omega),\omega) - b
[t]\ran}_{\cL_2^0}\\
\ns\ds \geq f(t,\bar x(t,\omega),u(t,\omega),\omega) - f [t]+2^{-i},
\ea
$$
and for a.e. $(t,\omega) \notin A_i$,
 $$
 \ba{ll}\ds
 {\lan P_{1}(t,\omega),a(t,\bar x(t,\omega),u(t,\omega),\omega) - a
[t]\ran}_H +{\lan Q_{1}(t,\omega),b(t,\bar x(t,\omega),u(t,\omega),\omega) - b
[t]\ran}_{\cL_2^0}\\
\ns\ds  =f(t,\bar x(t,\omega),u(t,\omega),\omega)  - f [t].
\ea
$$

Integrating the above two relations, we get a contradiction with   (\ref{integineq}). Consequently $\mu(A_i)=0$  for all $i.$\endpf

\begin{remark}
In  \cite{LZ3}  under stronger regularity assumptions and using the spike variation technique, the following maximality condition was proved :
\bel{eq1812061}
\ba{ll}\ds
{\langle P_1(t), a(t,\bar x(t),u) -   a[t])\rangle}_H  + {\langle Q_1(t), b(t,\bar x(t),u) -   b[t])\rangle}_{\cL_2^0} + f[t] - f(t,\bar x(t),u)  \\
\ns
\ds\q +\frac12 {\langle P_2(t)(b(t,\bar x(t),u) -   b[t]), b(t,\bar x(t),u) -   b[t]\rangle}_{\cL_2^0}\\
\ns
\ds  \leq 0, \q\;\; a.e. \; [0,T] \times \Omega,\; \forall u \in U,
\ea \ee
where $(P_1,Q_1)$  is given by Lemma \ref{vw-th1}, and $P_2(\cd)$ is provided by Theorem \ref{OP-th2}.

We claim that the inequality (\ref{eq1812061}) implies the conclusion of Theorem  \ref{PMPconvex}  whenever the sets $F(t,x,\omega)$  defined above are convex.   Indeed, fix $(t,\omega)$ such that the inequality (\ref{eq1812061}) is verified and let $u \in U$.  Then for all small $\varepsilon >0$, the vector
$$(a[t],b[t], f[t]) + \varepsilon (a(t,\bar x(t,\omega),u,\omega)-a[t],b(t,\bar x(t,\omega),u,\omega)-b[t], f(t,\bar x(t,\omega),u,\omega)-f[t] )  $$
is an element of $F(t,\bar x(t),\omega)$.
Consequently,  there exist $u_1 \in U$ and $ r \geq 0$ such that
$$
\ba{ll}
\ds a(t,\bar x(t,\omega),u_1,\omega)= a[t] + \varepsilon (a(t,\bar x(t,\omega),u,\omega)-a[t]),\\
\ns\ds
b(t,\bar x(t,\omega),u_1,\omega)= b[t] + \varepsilon (b(t,\bar x(t,\omega),u,\omega)-b[t]),\\
\ns\ds f(t,\bar x(t,\omega),u_1,\omega)+r = f[t]+ \varepsilon (f(t,\bar x(t,\omega),u,\omega)-f[t]).
\ea
$$
Hence, by the maximality condition and noting $r \geq 0$ we obtain that
$$
\ba{ll}\ds
{\langle P_1(t,\omega), a(t,\bar x(t,\omega),u_1,\omega) -   a[t]\rangle }_H + {\langle Q_1(t,\omega), b(t,\bar x(t,\omega),u_1,\omega) -   b[t]\rangle}_{\cL_2^0} + f[t] \\
 \ns\ds\q- f(t,\bar x(t,\omega),u_1,\omega) -r +\frac12 {\langle P_2(t,\omega)(b(t,\bar x(t,\omega),u_1,\omega) -   b[t]), b(t,\bar x(t,\omega),u_1,\omega) -   b[t]\rangle}_{\cL_2^0} \\\ns\ds \leq 0 .
 \ea
$$
This implies that
$$
\ba{ll}\ds
\varepsilon {\langle P_1(t,\omega), a(t,\bar x(t,\omega),u,\omega)-a[t]\rangle}_H  + \varepsilon {\langle Q_1(t,\omega), b(t,\bar x(t,\omega),u,\omega)-b[t]\rangle}_{\cL_2^0}\\
 \ns\ds\q
 - \varepsilon (f(t,\bar x(t,\omega),u,\omega)-f[t]) +\frac{\varepsilon^2}{2} {\langle P_2(t,\omega)(b(t,\bar x(t,\omega),u,\omega)-b[t]), b(t,\bar x(t,\omega),u,\omega)-b[t]\rangle}_{\cL_2^0}\\
 \ns\ds  \leq 0 .
 \ea$$
Dividing both sides of the above inequality by $\varepsilon$  and taking the limit we end the proof of our claim.
\end{remark}

\section{Second order necessary optimality condition}\label{s4}

In this section, we investigate the second order necessary conditions for the local minimizers $(\bx,\bu)$ of (\ref{jk2}).
In addition to the assumptions (A1)--(A3), we suppose that

\ms

\no{\bf (A6)} {\em The functions $a$, $b$, $f$ and $g$ satisfy the following:}
  \begin{enumerate}
      \item[i)] {\em For a.e. $(t,\omega)\in [0, T]\times\Omega$, the functions
             $a(t, \cdot, \cdot,\omega):\ H\times \wt H\to H$ and $b(t, \cdot, \cdot,\omega):\ H\times \wt H\to \cL_2^0$
             are twice differentiable
             and
             $$(x,u)\mapsto (a_{(x,u)^2}(t,x,u,\omega),b_{(x,u)^2}(t,x,u,\omega))$$
             is uniformly continuous in $x\in H$ and $u\in \wt H$,
             and,}
             $$|a_{(x,u)^2}(t,x,u,\omega)|_{\cL((H,\wt
H) \times (H,\wt
H);H)} +  |b_{(x,u)^2}(t,x,u,\omega)|_{\cL((H,\wt
H) \times (H,\wt
H);\cL_2^0)}\le L,\quad\forall\; (x,u)\in H\times \wt H;$$

       \item[ii)] {\em For a.e. $(t,\omega)\in [0, T]\times\Omega$, the functions $f(t, \cdot, \cdot,\omega):\ H\times \wt H\to \mr$ and $g(\cdot,\omega):\ H \to \mr$ are  twice continuously differentiable, and for any $x,\ \tilde{x}\in H$ and $u,\ \tilde{u}\in \wt H$,}
             $$
             \left\{
             \begin{array}{l}
              |f_{(x,u)^2}(t,x,u,\omega)|_{\cL((H,\wt
H) \times (H,\wt
H); \dbR)}\le L,\\
              |f_{(x,u)^2}(t,x,u,\omega)-f_{(x,u)^2}(t,\tilde{x},\tilde{u},\omega)|_{\cL((H,\wt
H) \times (H,\wt
H); \dbR)}\le L(|x-\tilde{x}|_H+|u-\tilde{u}|_{\wt H}),\\
              |g_{xx}(x,\omega)|_{\cL(H\times H;\dbR)} \le L, \ |g_{xx}(x,\omega)-g_{xx}(\tilde{x},\omega)|_{\cL(H\times H;\dbR)}\le L|x-\tilde{x}|_H.
             \end{array}\right.
             $$
  \end{enumerate}

For $\varphi=a,\; b, \;  f$, denote
\bel{181209e1}
\varphi_{xx}[t]=\varphi_{xx}(t,\bar{x}(t),\bar{u}(t)),\quad
\varphi_{xu}[t]=\varphi_{xu}(t,\bar{x}(t),\bar{u}(t)),\quad
\varphi_{uu}[t]=\varphi_{uu}(t,\bar{x}(t),\bar{u}(t)).
\ee

Let  $\bu, v, h,h_{\eps}\in L^{2\beta}_{\mmf}(\Omega; L^{4}(0,T;\wt H))$ ($\beta \geq 1$) be such that
\bel{191210e1}
|h_{\eps}-h|_{L^{2\beta}_{\mmf}(\Omega;L^{4}(0,T;\wt H))}\to0,\q  \hbox{ as }\eps\to 0^+.
 \ee
Set
 \bel{191210e2}
 u^{\eps}:=\bu+\eps v+\eps^2 h_{\eps}.
 \ee
Denote by $x^{\eps}$ the solution  of (\ref{ch-10-fsystem1})
corresponding to the control $u^{\eps}$ and the initial datum
$x_{0}$. Put
 \bel{191210e3}
 \delta x^{\eps}=x^{\eps}-\bar{x}, \qquad\delta u^{\eps}=\eps v+\eps^2 h_{\eps}.
 \ee

Similarly to \cite{Hoehener12}, we introduce the following second-order variational equation\footnote{Recall that, for any $C^2$-function $F (\cd):  X \to Y$ and $x_0 \in X$, $F_{xx}(x_0)\in\cL(X \times X; Y)$.  This means that, for any $x_1,x_2 \in X$, $F_{xx}(x_0)(x_1,x_2) \in Y$. Hence, by (\ref{181209e1}), $a_{xx}[t]\big(y_{1}(t),y_{1}(t)\big)$ (in (\ref{second order vari equ}))  stands for $a_{xx}(t,\bar{x}(t),\bar{u}(t))\big(y_{1}(t),y_{1}(t)\big)$.   One has a similar meaning for
$a_{uu}[t]\big(v(t),v(t)\big)$ and so on.}:
\begin{equation}\label{second order vari equ}
\quad\left\{
\begin{array}{l}
dy_{2}(t)= \Big(Ay_{2}(t) +a_{x}[t]y_{2}(t) + 2a_{u}[t]h(t) +a_{xx}[t]\big(y_{1}(t),y_{1}(t)\big)+2a_{xu}[t]\big(y_{1}(t),v(t)\big)\\
\qquad\qquad+a_{uu}[t]\big(v(t),v(t)\big)\Big)dt
+\Big(b_{x}[t]y_{2}(t)+2b_{u}[t]h(t) +b_{xx}[t]\big(y_{1}(t),y_{1}(t)\big)\\
\qquad\qquad+2b_{xu}[t]\big(y_{1}(t),v(t)\big)
+b_{uu}[t]\big(v(t),v(t)\big)\Big)dW(t),\qquad t\in[0,T],\\
y_{2}(0)=0,
\end{array}\right.
\end{equation}
where $y_{1}$ is the solution to the first variational equation (\ref{first vari equ1}) (for $v(\cdot)$ as above). The adjoint equation for  (\ref{second order vari equ}) is given by (\ref{second ajoint equ}) (See Theorem \ref{OP-th2} for its well-posedness in the sense of relaxed transposition solution).

Similarly to \cite[Lemma 4.1]{FZZ}, we have the following estimates for solutions to (\ref{second order vari equ}).
 \begin{lemma}\label{estimate two of varie qu}
Assume (A1)--(A3),  (A6)  and let $\beta \geq 1$. Then, for  $\bu, v, h,h_{\eps}\in L^{2\beta}_{\mmf}(\Omega; L^{4}(0,T;\wt H))$ so that (\ref{191210e1}) holds and for $ \delta x^{\eps}$ given by (\ref{191210e3}), we have
\begin{equation*}
\|y_{2}\|_{L_{\mmf}^{\infty}(0,T; L^{\beta}(\Omega; H))}\le C\left(\|v\|_{L^{2\beta}_{\mmf}(\Omega; L^{4}(0,T;\wt H))}+\|h\|_{L^{\beta}_{\mmf}(\Omega; L^{2}(0,T;\wt H))}\right).
\end{equation*}
Furthermore,
\begin{equation}\label{r2 to 0}
\|r_{2}^{\eps}\|_{L_{\mmf}^{\infty}(0,T; L^{\beta}(\Omega; H))}\to 0,\quad \hbox{ as }\eps\to 0^+,
\end{equation}
where,
$$r_{2}^{\eps}(t,\omega):=\frac{\delta x^{\eps }(t,\omega)-\eps  y_{1}(t,\omega)}{\eps ^2} -\frac{1}{2}y_{2}(t,\omega).$$
\end{lemma}

\def\mS{\mathrm{S}}
\def\ms{\mathbb{S}}

To simplify the notation, we define\footnote{Note that the definition of $\ms(t,x,u,y_{1},z_{1},y_{2},\omega)$ in (\ref{S}) is different from that in \cite[the equality (4.4), p. 3708]{FZZ}. The main reason for this is due to the fact that the characterization of $Q(\cd)$ in Theorem \ref{OP-th2} is much weaker than the one  in the finite dimensions.} (Recall (\ref{Hamiltonianconvex}) for the definition of $\mathbb{H}(\cd)$):
\begin{equation}\label{S}
\ba{ll}
\ds \ms(t,x,u,y_{1},z_{1},y_{2},\omega)
:=&\ds\mathbb{H}_{xu}(t,x,u,y_{1},z_{1},\omega)+a_{u}(t,x,u,\omega)^*y_{2}\\[3mm]
&\ds+b_{u}(t,x,u,\omega)^*
y_{2}b_{x}(t,x,u,\omega),
\ea
\end{equation}
where $(t,x,u,y_{1},z_{1},y_{2},\omega)\in [0,T]\times H\times \wt H\times H\times \cL_2^0\times \cL(H)\times\Omega$,
and denote
\begin{equation}\label{S(t)}
\ms[t]= \ms(t,\bar{x}(t),\bar{u}(t),P_{1}(t),Q_{1}(t),P_{2}(t)),
\quad t \in [0,T].
\end{equation}

Let $\bu\in \mmu\cap L_{\mmf}^{4}(0,T;\wt H)$. Define the critical set
$$\ups :=\left\{v \in L_{\mmf}^{2}(0,T;\wt H)) \ \Big|\
\lan \mathbb{H}_{u}(t,\omega),v(t,\omega){\ran}_{\wt H}=0 \;\;\mbox{\rm a.e.} \;\; t \in [0,T],   \;\;\mbox{\rm a.s.}\right\},$$
and the set of admissible second order variations by
\begin{eqnarray*}
& &\ma:=\Big\{(v, h)\in L_{\mmf}^{4}(0,T;\wt H))\times L_{\mmf}^{4}(0,T;\wt H))\ ~\Big|~\
\\
& &\qquad\qquad\qquad\qquad
 h(t,\omega)\in \TT_{U}(\bar u (t,\omega), v(t,\omega)), \
\mbox{a.e.}\ t\in[0,T],\ \mbox{a.s. }
\Big\}.
\end{eqnarray*}

We have the following result.

\begin{theorem}\label{TH second order integral condition}
Let (A1)--(A3) and (A6) hold, $(\bx,\bu)$ be an optimal pair for for Problem (OP) and $\bu\in L_{\mmf}^{4}(0,T;\wt H))$.   Then,
for all $(v,h)\in \ma$ with $v\in \ups$, it holds that
\begin{equation}\label{second order integral condition}
\begin{array}{ll}
\ds\dbE\int_{0}^{T}\Big[2\lan \mathbb{H}_{u}[s],h(s){\ran}_{\wt H}+2\lan \ms[s]y_{1}(s),v(s){\ran}_{\wt H}
+\lan \mathbb{H}_{uu}[s]v(s),v(s){\ran}_{\wt H}\\
\ns\ds \qq+\big\langle  P_2(s)b_{u}[s]v(s), b_{u}[s]v(s)\big\rangle_{\cL_2^0}+\big\langle b_{u}[s]v(s), \widehat
Q_2^{(0)}(0,a_{u}v,b_{u}v)(s)\big\rangle_{\cL_2^0}\\
\ns\ds \qq+
\big\langle
Q_2^{(0)}(0,a_{u}v,b_{u}v)(s), b_{u}[s]v(s)
\big\rangle_{\cL_2^0}\Big]ds\\
\ns\ds\le 0.
\end{array}
\end{equation}
\end{theorem}

{\it Proof.} Since $\bu\in L_{\mmf}^{4}(0,T;\wt H)$, for any $(v,h)\in \ma$ with $v\in \ups$, similarly to the proof of \cite[Theorem 4.1]{FZZ}, without loss of generality, we may choose $h_{\eps}\in L_{\mmf}^4(0,T;\wt H))$ so that (\ref{191210e1}) holds for $\beta=2$. Let $ u^{\eps}$,
$x^{\eps}$, $\delta x^{\eps}$ and $\delta u^{\eps}$ be  as  in (\ref{191210e2})--(\ref{191210e3}). Denote
$$\tilde{f}_{xx}^{\eps}(t):=\int_{0}^{1}(1-\theta)f_{xx}(t,\bx(t)
+ \theta\delta x^{\eps}(t),\bu(t)+\theta\delta
u^{\eps}(t))d\theta.
$$
Mappings  $\tilde{f}_{xu}^{\eps}(t)$,
$\tilde{f}_{uu}^{\eps}(t)$ and $\tilde{g}_{xx}^{\eps}(T)$ are
defined in a similar way.

Expanding the cost functional $J$ at $\bu$, we get
\begin{eqnarray*}
& & \frac{J(u^{\eps})-J(\bu)}{\eps^2}\\
&&=\frac{1}{\eps^2}\dbE\int_{0}^{T}\Big(\lan f_x[t],\delta x^{\eps}(t){\ran}_{H}+
\lan f_{u}[t],\delta u^{\eps}(t){\ran}_{\wt H}
+\lan\tilde{f}_{xx}^{\eps}(t)\delta x^{\eps}(t),\delta x^{\eps}(t){\ran}_{H}
\\
& &\q
+2\lan\tilde{f}_{xu}^{\eps}(t)\delta x^{\eps}(t),\delta u^{\eps}(t){\ran}_{\wt H}
+\lan\tilde{f}_{uu}^{\eps}(t)\delta u^{\eps}(t),\delta u^{\eps}(t){\ran}_{\wt H}\Big)dt\\
& &
\q+\frac{1}{\eps^2}\dbE \Big(\lan g_{x}(\bx(T)),\delta x^{\eps}(T){\ran}_{H}+\lan \tilde{g}_{xx}^{\eps}(\bx(T))\delta x^{\eps}(T),\delta x^{\eps}(T){\ran}_{H}
\Big)\\
&&= \dbE\int_{0}^{T}\Big[\frac{1}{\eps}\lan f_{x}[t],y_{1}(t){\ran}_{H}
+\frac{1}{2}\lan f_{x}[t],y_{2}(t){\ran}_{H}
+\frac{1}{\eps}\lan f_{u}[t],v(t){\ran}_{\wt H}+\lan f_{u}[t],h(t){\ran}_{\wt H}\\
& &\q+\frac{1}{2}\Big(\lan f_{xx}[t]y_{1}(t),y_{1}(t){\ran}_{H}
+2\lan f_{xu}[t]y_{1}(t),v(t){\ran}_{\wt H}
+\lan f_{uu}[t]v(t),v(t){\ran}_{\wt H}\Big)\Big]dt\\
& &\q
+\dbE \Big(\frac{1}{\eps} \lan g_{x}(\bar{x}(T)),y_{1}(T){\ran}_{H}
+\frac{1}{2} \lan g_{x}(\bx(T)),y_{2}(T){\ran}_{H} \\
& &\q+\frac{1}{2}\lan g_{xx}(\bx(T))y_{1}(T),y_{1}(T){\ran}_{H}\Big)
+ \rho_{2}^{\eps},
\end{eqnarray*}
where
\begin{eqnarray*}
\rho_{2}^{\eps}
&=&\dbE \int_{0}^{T}\Big( \lan f_{x}[t],r_{2}^{\eps}(t){\ran}_{H}+
\lan f_{u}[t], h_{\eps}(t)-h(t){\ran}_{\wt H}\Big)dt+\dbE\lan g_{x}(\bx(T)),r_{2}^{\eps}(T){\ran}_{H}\\
&&
+\dbE \int_{0}^{T}\Big[\Big(\left\langle \tilde{f}_{xx}^{\eps}(t)\frac{\delta x^{\eps}(t)}{\eps},\frac{\delta x^{\eps}(t)}{\eps}\right\rangle_{H}-\frac{1}{2}\lan f_{xx}[t]y_{1}(t),y_{1}(t){\ran}_{H}\Big)\\
&&\qquad\quad
+\Big(2\left\langle \tilde{f}_{xu}^{\eps}(t)\frac{\delta x^{\eps}(t)}{\eps},\frac{\delta u^{\eps}(t)}{\eps}\right\rangle_{\wt H}-\lan f_{xu}[t]y_{1}(t),v(t){\ran}_{\wt H}\Big)\\
&&\qquad\quad
+\Big(\left\langle\tilde{f}_{uu}^{\eps}(t)\frac{\delta u^{\eps}(t)}{\eps},\frac{\delta u^{\eps}(t)}{\eps}\right\rangle_{\wt H}-\frac{1}{2}\lan f_{uu}[t]v(t),v(t){\ran}_{\wt H}\Big)\Big]dt\\
&&+\dbE \Big(\left\langle\tilde{g}_{xx}^{\eps}(\bx(T))\frac{\delta x^{\eps}(T)}{\eps},\frac{\delta x^{\eps}(T)}{\eps}\right\rangle_{H}-\frac{1}{2}\lan g_{xx}(\bx(T))y_{1}(T),y_{1}(T){\ran}_{H}\Big).
\end{eqnarray*}

As in Lemma \ref{estimate two of varie qu}, we  find that $\lim_{\eps\to 0^+}\rho_{2}^{\eps}=0$. On the other hand,
by calculations done in  \eqref{first order taylor exp}--\eqref{first order vari inequ} and the definition of $\mathbb{H}(\cd)$) in (\ref{Hamiltonianconvex}),  and recalling that $v\in \ups$, we have
\begin{eqnarray*}
& &\frac{1}{\eps}\dbE\int_{0}^{T}\Big(\lan f_{x}[t],y_{1}(t){\ran}_{H}
+\lan f_{u}[t],v(t){\ran}_{\wt H}\Big)dt
+\frac{1}{\eps}\dbE\lan g_{x}(\bar{x}(T)),y_{1}(T){\ran}_{H}\nonumber\\
&&=-\frac{1}{\eps}\dbE\int_{0}^{T}\lan \mathbb{H}_{u}[t],v(t){\ran}_{H}dt=0.
\end{eqnarray*}
Therefore,
\begin{equation}\label{taylorexpconvex}
\begin{array}{ll}
\displaystyle 0\le\lim_{\eps\to 0^+} \frac{J(u^{\eps}(\cdot))-J(\bu(\cdot))}{\eps^2}\\[3mm]\displaystyle
\quad= \dbE\int_{0}^{T}\Big[\frac{1}{2}\lan f_{x}[t],y_{2}(t){\ran}_{H}
+\lan f_{u}[t],h(t){\ran}_{\wt H}\\[3mm]\displaystyle\qquad
+\frac{1}{2}\Big(\lan f_{xx}[t]y_{1}(t),y_{1}(t){\ran}_{H}
+2\lan f_{xu}[t]y_{1}(t),v(t){\ran}_{\wt H}
+\lan f_{uu}[t]v(t),v(t){\ran}_{\wt H}\Big)\Big]dt\\[3mm]\displaystyle\qquad
+\frac{1}{2}\dbE\Big(\lan g_{x}(\bx(T)),y_{2}(T){\ran}_{H}
+\lan g_{xx}(\bx(T))y_{1}(T), y_{1}(T){\ran}_{H}\Big).
\end{array}
\end{equation}

From \eqref{eq def solzz} and \eqref{second order vari equ}, it follows that
\begin{equation}\label{hxty2}
\begin{array}{ll}
\dbE \lan g_{x}(\bar{x}(T)),y_{2}(T){\ran}_{H}\\[3mm]\displaystyle
=\dbE\int_0^T \big\langle a_x[s]^* P_1(s) +b_x[s]^* Q_1(s)- f_x[s],y_2(s)\big\rangle_Hds\\[3mm]\displaystyle
\quad -\dbE\int_0^T \big\langle
P_1(s),a_{x}[s]y_{2}(s) + 2a_{u}[s]h(s) +a_{xx}[s]\big(y_{1}(s),y_{1}(s)\big)\\
\quad  +2a_{xu}[s]\big(y_{1}(s),v(s)\big)+a_{uu}[s]\big(v(s),v(s)\big)\big\rangle_H ds \\[3mm]\displaystyle
\quad  - \dbE\int_0^T
\big\langle Q_1(s),b_{x}[s]y_{2}(s)+2b_{u}[s]h(s) +b_{xx}[s]\big(y_{1}(s),y_{1}(s)\big)\\[3mm]\displaystyle
\quad  +2b_{xu}[s]\big(y_{1}(s),v(s)\big)
+b_{uu}[s]\big(v(s),v(s)\big)\big\rangle_{\cL_2^0} ds\\[3mm]\displaystyle
=-\dbE\int_0^T \Big[\big\langle  f_x[s],y_2(s)\big\rangle_H+ \big\langle
P_1(s), 2a_{u}[s]h(s) +a_{xx}[s]\big(y_{1}(s),y_{1}(s)\big)\\[3mm]\displaystyle
\quad  +2a_{xu}[s]\big(y_{1}(s),v(s)\big)+a_{uu}[s]\big(v(s),v(s)\big)\big\rangle_H   +
\big\langle Q_1(s),2b_{u}[s]h(s) +b_{xx}[s]\big(y_{1}(s),y_{1}(s)\big)\\[3mm]\displaystyle
\quad  +2b_{xu}[s]\big(y_{1}(s),v(s)\big)
+b_{uu}[s]\big(v(s),v(s)\big)\big\rangle_{\cL_2^0}\Big] ds.
\end{array}\end{equation}
On the other hand, by \eqref{6.18eq1} and \eqref{first vari equ1}, we find that
\begin{equation}\label{hxxty12}
\begin{array}{ll}
\ds\mE\big\langle g_{xx}(\bar{x}(T)) y_1(T), y_1(T)
\big\rangle_{H} \\
\ns\ds =\mE \int_0^T \big\langle
\mathbb{H}_{xx}[s] y_1(s), y_1(s) \big\rangle_{H}ds- \mE \int_0^T
\big\langle P_2(s)a_{u}[s]v(s), y_1(s)\big\rangle_{H}ds
\\
\ns\ds \q-\mE\int_0^T\big\langle
P_2(s)y_1(s),
a_{u}[s]v(s)\big\rangle_{H}ds -\mE\int_0^T
\big\langle P_2(s)b_x[s]y_1 (s), b_{u}[s]v(s)\big\rangle_{\cL_2^0}ds \\
\ns\ds \q- \mE
\int_0^T\big\langle  P_2(s)b_{u}[s]v(s), b_x [s]y_1 (s) + b_{u}[s]v(s)\big\rangle_{\cL_2^0}ds\\
\ns\ds \q- \mE\int_0^T
\big\langle b_{u}[s]v(s), \widehat
Q_2^{(0)}(0,a_{u}v,b_{u}v)(s)\big\rangle_{\cL_2^0}ds-
\mE\int_0^T\big\langle
Q_2^{(0)}(0,a_{u}v,b_{u}v)(s), b_{u}[s]v(s)
\big\rangle_{\cL_2^0}ds\\
\ns\ds =\mE \int_0^T\Big[ \big\langle
\mathbb{H}_{xx}[s] y_1(s), y_1(s) \big\rangle_{H}- 2
\big\langle P_2(s)a_{u}[s]v(s), y_1(s)\big\rangle_{H}-2
\big\langle P_2(s)b_x[s]y_1 (s), b_{u}[s]v(s)\big\rangle_{\cL_2^0}\\
\ns\ds \q- \big\langle  P_2(s)b_{u}[s]v(s), b_{u}[s]v(s)\big\rangle_{\cL_2^0}- \big\langle b_{u}[s]v(s), \widehat
Q_2^{(0)}(0,a_{u}v,b_{u}v)(s)\big\rangle_{\cL_2^0}\\
\ns\ds \q-
\big\langle
Q_2^{(0)}(0,a_{u}v,b_{u}v)(s), b_{u}[s]v(s)
\big\rangle_{\cL_2^0}\Big]ds.
\end{array}
\end{equation}
Substituting (\ref{hxty2}) and (\ref{hxxty12}) into
(\ref{taylorexpconvex})  yields
$$
\begin{array}{ll}
\displaystyle 0\ge-\dbE\int_{0}^{T}\Big[\lan f_{u}[s],h(s){\ran}_{\wt H}
+\lan f_{xu}[s]y_{1}(s),v(s){\ran}_{\wt H}
+\frac{1}{2}\lan f_{uu}[s]v(s),v(s){\ran}_{\wt H}\Big]ds\\[3mm]\displaystyle\qquad
+\dbE\int_0^T \Big[\big\langle
P_1(s), a_{u}[s]h(s) +a_{xu}[s]\big(y_{1}(s),v(s)\big)+\frac{1}{2}a_{uu}[s]\big(v(s),v(s)\big)\big\rangle_H  \\[3mm]\displaystyle
\qquad +
\big\langle Q_1(s),b_{u}[s]h(s) +b_{xu}[s]\big(y_{1}(s),v(s)\big)
+\frac{1}{2}b_{uu}[s]\big(v(s),v(s)\big)\big\rangle_{\cL_2^0}\Big] ds\\[3mm]\displaystyle\qquad
+\frac{1}{2}\mE \int_0^T\Big[  2
\big\langle P_2(s)a_{u}[s]v(s), y_1(s)\big\rangle_{H}+2
\big\langle P_2(s)b_x[s]y_1 (s), b_{u}[s]v(s)\big\rangle_{\cL_2^0}\\
\ns\ds \qq+\big\langle  P_2(s)b_{u}[s]v(s), b_{u}[s]v(s)\big\rangle_{\cL_2^0}+\big\langle b_{u}[s]v(s), \widehat
Q_2^{(0)}(0,a_{u}v,b_{u}v)(s)\big\rangle_{\cL_2^0}\\
\ns\ds \qq+
\big\langle
Q_2^{(0)}(0,a_{u}v,b_{u}v)(s), b_{u}[s]v(s)
\big\rangle_{\cL_2^0}\Big]ds\\[3mm]\displaystyle
=\frac{1}{2}\dbE\int_{0}^{T}\Big[2\lan \mathbb{H}_{u}[s],h(s){\ran}_{\wt H}+2\lan \ms[s]y_{1}(s),v(s){\ran}_{\wt H}
+\lan \mathbb{H}_{uu}[s]v(s),v(s){\ran}_{\wt H}\\
\ns\ds \qq+\big\langle  P_2(s)b_{u}[s]v(s), b_{u}[s]v(s)\big\rangle_{\cL_2^0}+\big\langle b_{u}[s]v(s), \widehat
Q_2^{(0)}(0,a_{u}v,b_{u}v)(s)\big\rangle_{\cL_2^0}\\
\ns\ds \qq+
\big\langle
Q_2^{(0)}(0,a_{u}v,b_{u}v)(s), b_{u}[s]v(s)
\big\rangle_{\cL_2^0}\Big]ds.
\end{array}
$$
Then, we obtain the desired inequality
(\ref{second order integral condition}). This completes the proof of Theorem \ref{TH second
order integral condition}.
\endpf

\end{document}